\documentclass[a4paper,oneside,notitlepage]{article}%
\usepackage{amsmath}
\usepackage{amsfonts}
\usepackage{amssymb}
\usepackage{graphicx}
\usepackage{hyperref}%
\pdfoutput=1
\providecommand{\U}[1]{\protect\rule{.1in}{.1in}}
\newtheorem{theorem}{Theorem}

\newtheorem{lemma}{Lemma}

\newtheorem{proposition}{Proposition}

\newenvironment{proof}[1][Proof]{\noindent\textbf{#1.} }{\ \rule{0.5em}{0.5em}}
\begin{document}

\title{A space-time stochastic model and a covariance function for the stationary
spatio-temporal random process and spatio-temporal prediction (kriging)}
\author{T. Subba Rao\\University of Manchester, U.K. and \\C. R. Rao AIMSCS, University of Hyderabad Campus, India,\\tata.subbarao$@$gmail.com
\and Gyorgy Terdik\\University of Debrecen, Hungary, \\Terdik.Gyorgy$@$inf.unideb.hu }
\date{}
\maketitle

\begin{abstract}
Consider a stationary spatio-temporal random process $\{Y_{t}(\mathbf{s}%
)\mathbf{|s}\in\mathbb{R}^{d},$ $\ t\in\mathbb{Z\}}$ and let $\{Y_{t}%
(\mathbf{s}_{i})\mathbf{|}i=1,2,...m;\ t=1,\ldots,n$ $\mathbb{\}}$ be a sample
from the random process$\ \{Y_{t}\left(  \mathbf{s}\right)  \mathbf{\}}$. Our
object here is to estimate $\{Y_{t}(\mathbf{s}_{o)}\}$ for all $t$ at the
location $\mathbf{s}_{o}$ given the sample $\{Y_{t}(\mathbf{s}_{i}%
)|\ i=1,2,...m;\ t=1,\ldots,n\mathbf{\}}$ using the frequency domain approach.
To obtain an estimator, we define a sequence of Discrete Fourier transforms
$\{J_{\mathbf{s}}(\omega_{i})\}$ at \ the Fourier frequencies using the time
series observed at the locations $\mathbf{s}_{i},$ $\left(
i=1,2,3...m\right)  $, and use these complex random variables as our
observations from the complex valued random process $\{J_{\mathbf{s}}%
(\omega)\}$. The \ Fourier transforms are functions of the spatial coordinates
only.\ Assuming the complex valued process satisfies a complex stochastic
partial differential equation (CSPDE) of the Laplacian type, and using the
properties of Fourier transforms of stationary processes, we obtain an
expression for the spatio-temporal covariance function and \ for the spectral
density function. The covariance function of the Discrete Fourier transforms
at two distinct locations has been shown to be a function of the Euclidean
\ distance and the temporal frequency \ and its second order spectrum
corresponds to non separable class of random processes. We show further that
the model defined here includes as special cases the spatio-temporal models
defined by \cite{jones1997models}, \cite{lindgren2011explicit} and
\cite{sigrist2015stochastic}. The estimation \ of the parameters of the
spatio-temporal covariance function \ has also been considered. The method of
estimation is based on Frequency Variogram approach recently introduced and it
does not involve inversion of large dimensional matrices and found to be
robust against departure from Gaussianity.\ Since the methods developed here
are based on the Discrete Fourier Transforms, and \ these transforms can be
evaluated using Fast Fourier Transform algorithms, the estimation methods are
quick and efficient compared to the classical approaches based on Variogram
method and\ the likelihood approaches. The methods are illustrated with a real
example. The data considered is Air Pollution data (Particulate Matter
PM$_{2.5}$) recorded at 15 locations in New York city observed\ over a period
of 9 months and the data at the location 11 has large number of missing
values. The data at this location is estimated and prediction intervals are
given and the estimates are compared with the observed data when the data is
available. \bigskip


\end{abstract}

\section{Introduction and Summary}

In recent years it has become necessary to develop statistical methods for the
analysis of data coming from diverse areas such as, environment, marine
biology, agriculture, finance etc. The data which comes from these areas, are
usually, functions of both space and time. Any statistical method developed
must take into account \ both spatial dependence, temporal dependence and any
interaction between space and time. There is a vast literature on statistical
analysis of stationary spatial data (for example refer to the books of
\cite{Cressie1993}, \cite{stein1999interpolation}) but not \ to the same
extent in the case of stationary spatio-temporal data. The inclusion of an
extra temporal dimension, which cannot be imbedded into spatial dimension
gives raise to many problems. One such problem is finding a suitable
covariance function \ which is positive definite and depends on spatial lag
difference and temporal lag. In recent years several authors (see
\cite{cressie1999classes}, \cite{gneiting2002nonseparable},
\cite{diggle2007model}, \cite{stein2005space}, \cite{craigmile2011space},
\cite{jones1997models}, \cite{ma2002spatio}, \cite{ma2003spatio},
\cite{lindgren2011explicit}, and \cite{sigrist2015stochastic}) have proposed
various types of covariance functions, and majority of them are Matern type of
functions. \cite{jones1997models}, \cite{lindgren2011explicit} and
\cite{sigrist2015stochastic} have considered transport-diffusion Stochastic
Partial Differential (SPDE) equations for modelling of stationary
spatio-temporal random process. The models \ defined by these authors are
stochastic versions of the \ classical Heat equation, which is in fact a
dynamic form of the Laplace equation, and the temporal correlation is
explained by inclusion of a first order time derivative. In a way this is
equivalent to assuming that always the temporal correlation in the data can be
explained by a first order Autoregressive model(AR(1)). This assumption some
times can be very unrealistic. \cite{lindgren2011explicit} considered the
approximation of Gaussian Field (GF) models by the Gaussian \ Markov Random
Field (GMRF) models and considered modelling the data by GMRF models.
\cite{sigrist2015stochastic} have approximated the solution of the SPDE models
by a linear combination of deterministic spatial functions (Fourier functions
\ in terms of spatial wave numbers) with random coefficients that evolve
dynamically and requires discretization \ in time for application to
\ discrete time series.

The prediction of the data at a known location (usually known as Kriging) is
another important problem. The object here is to predict the data at a known
location where the time series data is not observed. It is well known that the
best linear minimum mean square predictor depends on the unknown \ covariance
function of the process . The elements of \ the variance covariance matrix
cannot be estimated from the data as no data is available at the location
where we are estimating the data. Therefore, there is a need to find a
suitable parametric covariance function. Our main object in this paper is to
explore and investigate possible ways of using the Discrete Fourier Transforms
of spatio-temporal data in modelling and prediction(for further applications
of the DFT to spatial data analysis, we refer to \cite{raoGyT2012HB}. We use
many well known asymptotic sampling properties of the Discrete Fourier
Transforms of stationary time series\ in our derivations.

The novel features of our present paper are as follows. We consider the
discrete Fourier transforms of the given\ spatio-temporal data and treat these
complex Gaussian variables as our data. We define the complex Gaussian process
through a stochastic partial differential equation in spatial coordinates
only. By defining in this way and operating on the Fourier transforms, the
dependency of the operator on time is removed and hence no discretization in
time is necessary.\ By analyzing the Complex Gaussian processes and using the
defined Complex Stochastic Partial Differential Equations (CSPDE), we are
reducing the number of computations required for solving the equations, the
estimation of the parameters\ of the covariance function and for forecasting
too. Under the assumption of isotropy, we obtain an expression for the
covariance of the Discrete Fourier Transforms at different locations in terms
of the Euclidean distance between the locations, and it is shown to be in
terms of the modified Bessel function (Matern Class).The covariance function
obtained here is a function of the Euclidean spatial distance, and the
temporal frequency spectrum characterising the temporal dependence in the
series In this way the expression for the covariance function given here is
fundamentally different from other covariance functions defined and used by
other authors. We further show that the \ second order spectral density
function of the spatio-temporal random process defined here through the above
operator on the discrete Fourier transforms includes all the second order
spectra of the process so far defined \ through the stochastic versions of the
Laplace equation such as, Heat equation (transport-diffusion equation of
\cite{jones1997models}, \cite{lindgren2011explicit} and
\cite{sigrist2015stochastic}), Wave equation and Helmholtz equation as special
cases. The second order spectral density function obtained here belongs to
nonseparable random process.

We summarize the contents of each section. In section 2, the notation and the
spectral representation of the spatio-temporal random process \ are
introduced. The properties of \ the Discrete Fourier transforms of second
order stationary processes are discussed in the appendix. Expressions for the
spatio-temporal covariance and for the spectral density functions when the
Discrete Fourier Transforms satisfy a Complex Stochastic Partial Differential
equation (CSPDE) are obtained in section 3. We also show in this section that
the second order spectra of the processes satisfying the models defined by
\cite{jones1997models}, \ and \cite{sigrist2015stochastic} can be obtained as
special cases of the CSPDE model defined here. The main results related to the
process satisfying CSPDE are stated in theorems 1 and 2. The prediction of the
entire data at a known location given the data in the neighborhood using the
Discrete Fourier transforms is considered in section 4.\ The estimation of the
parameters, using Frequency Variogram method, \ of the spatio-temporal
covariance function is considered in section 5.\ In section 6, the analysis of
the Air Pollution data (Particulate Matter PM$_{2.5}$) collected at 15
locations in New York City is considered. The prediction of the entire data at
the location 11, where large number of observations are missing \ is also
considered in section 6 and the plots of the estimated data together with
prediction intervals are also given. The estimates of the observations
corresponding to PM$_{2.5}$ \ at time $t=22$ \ and at all the 625 locations
(the spatial coordinates are known) are calculated and the surface plot of the
estimated values is shown in Fig 4 of section 6.

\section{Notation and Preliminaries}

Let $Y_{t}\left(  \mathbf{s}\right)  $, where $\left\{  \mathbf{s}%
\in\mathbb{R}^{d},t\in\mathbb{Z}\right\}  $, denote the spatio-temporal random
process. We assume that the random process is spatially and temporally second
order stationary, i.e.
\begin{align*}
E\left[  Y_{t}\left(  \mathbf{s}\right)  \right]   &  =\mu,\\
Var\left[  Y_{t}\left(  \mathbf{s}\right)  \right]   &  =\sigma_{Y}^{2}%
<\infty,\\
Cov\left[  Y_{t}\left(  \mathbf{s}\right)  ,Y_{t+u}\left(  \mathbf{s}%
+\mathbf{h}\right)  \right]   &  =c\left(  \mathbf{h},u\right)  ,\quad
\mathbf{h}\in\mathbb{R}^{d},\,u\in\mathbb{Z}\mathbf{.}%
\end{align*}
We note that $c\left(  \mathbf{h},0\right)  $ and $c\left(  \mathbf{0}%
,u\right)  $ correspond to the purely spatial and purely temporal covariances
of the process respectively. A further common stronger assumption that is
often made is that the process \ is isotropic. The assumption of isotropy is a
stronger assumption. The process is said to be isotropic if
\[
c\left(  \mathbf{h},u\right)  =c\left(  \left\Vert \mathbf{h}\right\Vert
,u\right)  ,\quad\left\Vert \mathbf{h}\right\Vert \geq0,\,u\in\mathbb{Z},
\]
where $\left\Vert \mathbf{h}\right\Vert $ is the Euclidean distance. Without
loss of generality, we set $\mu$ equal to zero. As in the case of spatial
process, one can define the spatio-temporal variogram for $\left\{
Y_{t}\left(  \mathbf{s}\right)  \right\}  $ as
\begin{equation}
2\gamma\left(  \mathbf{h},u\right)  =Var\left[  Y_{t+u}\left(  \mathbf{s}%
+\mathbf{h}\right)  -Y_{t}\left(  \mathbf{s}\right)  \right]  .
\label{eq 1.1.}%
\end{equation}
If the random process $\left\{  Y_{t}\left(  \mathbf{s}\right)  \right\}  $ is
spatially and temporally stationary, then we can rewrite the above as
\begin{equation}
2\gamma\left(  \mathbf{h},u\right)  =2\left[  c\left(  \mathbf{0},0\right)
-c\left(  \mathbf{h},u\right)  \right]  , \label{eq 1.2.}%
\end{equation}
and for an isotropic process, $\gamma\left(  \mathbf{h},u\right)
=\gamma\left(  \left\Vert \mathbf{h}\right\Vert ,u\right)  $. We note\ that
$\gamma\left(  \mathbf{h},u\right)  $ is defined as the semi-variogram. \

In view of our assumption that the zero mean random process $\left\{
Y_{t}\left(  \mathbf{s}\right)  \right\}  $ is second order spatially and
temporally stationary, we have the spectral representation%
\begin{equation}
Y_{t}\left(  \mathbf{s}\right)  =\int\limits_{-\infty}^{\infty}\int
\limits_{-\pi}^{\pi}e^{i\left(  \mathbf{s}\cdot\underline{\lambda}%
+t\omega\right)  }dZ_{Y}\left(  \underline{\lambda},\omega\right)  ,
\label{eq 1.3.}%
\end{equation}
where $\mathbf{s}\cdot\underline{\lambda}=\sum\limits_{i=1}^{d}\mathbf{s}%
_{i}\lambda_{i}$ and $\int\limits_{-\infty}^{\infty}$ represents $d-$fold
multiple integral. We note that $Z_{Y}\left(  \underline{\lambda}%
,\omega\right)  $ is a zero mean complex valued random process with orthogonal
increments with
\begin{align}
E\left[  dZ_{Y}\left(  \underline{\lambda},\omega\right)  \right]   &
=0,\nonumber\\
E\left\vert dZ_{Y}\left(  \underline{\lambda},\omega\right)  \right\vert ^{2}
&  =dF_{Y}\left(  \underline{\lambda},\omega\right)  , \label{eq 1.4.}%
\end{align}
where $dF_{Y}\left(  \underline{\lambda},\omega\right)  $ is a spectral
measure. If we assume further that $dF\left(  \underline{\lambda}%
,\omega\right)  $ is absolutely continuous with respect to Lebesgue measure
$\ $then $dF\left(  \underline{\lambda},\omega\right)  =f\left(
\underline{\lambda},\omega\right)  \,d\underline{\lambda}\,d\omega$. Here
$f\left(  \underline{\lambda},\omega\right)  $ which is strictly positive and
real valued, is defined as the spatio-temporal spectral density function of
the random process $\left\{  Y_{t}\left(  \mathbf{s}\right)  \right\}  $, and
$-\infty<\lambda_{1},\ \lambda_{2},\ \ldots,\ \lambda_{d}<\infty,\ -\pi
\leq\omega\leq\pi$. In view of the orthogonality of the function $Z_{Y}\left(
\underline{\lambda},\omega\right)  $, we can show that the positive definite
covariance function $c\left(  \mathbf{h},u\right)  $ has the representation%
\begin{equation}
c\left(  \mathbf{h},u\right)  =\int\limits_{-\infty}^{\infty}\int
\limits_{-\pi}^{\pi}e^{i\left(  \mathbf{h}\cdot\underline{\lambda}%
+u\omega\right)  }f\left(  \underline{\lambda},\omega\right)  \,d\underline
{\lambda}\,d\omega, \label{eq 1.5.}%
\end{equation}
and by Fourier inversion, we have
\begin{equation}
f\left(  \underline{\lambda},\omega\right)  =\frac{1}{\left(  2\pi\right)
^{d+1}}\sum\limits_{u}\int\limits_{-\infty}^{\infty}e^{-i\left(
\mathbf{h\cdot}\underline{\lambda}+u\omega\right)  }c\left(  \mathbf{h}%
,u\right)  \,d\mathbf{h}, \label{eq 1.6.}%
\end{equation}
where $d\mathbf{h}=\prod\limits_{i=1}^{d}dh_{i}$. From \ equation
(\ref{eq 1.5.}) we obtain%
\[
c\left(  0,u\right)  =\int\limits_{-\infty}^{\infty}\int\limits_{-\pi}^{\pi
}e^{iu\omega}f\left(  \underline{\lambda},\omega\right)  \,d\underline
{\lambda}\,d\omega=\int\limits_{-\pi}^{\pi}e^{iu\omega}g_{0}(\omega)d\omega,
\]
where $g_{0}(\omega)=\int f\left(  \underline{\lambda},\omega\right)
\,d\underline{\lambda}$ is the temporal spectrum of the spatio-temporal random
process $\{Y_{t}\left(  \mathbf{s}\right)  \}$. In view of our assumption of
spatial stationarity $g_{0}(\omega)$ is same for all $\mathbf{s}$. It may be
pointed out here that a study of the properties of the second order temporal
spectrum in the spatio-temporal context can be of considerable interest in
several scientific fields, for example in neurosciences as shown by
\cite{ombao2008spatio}. \cite{ombao2008spatio} considered the estimation of
temporal spectrum at a given location assuming that spatio-temporal process is
slowly spatially changing using a methodology \ similar to
\cite{priestley1965evolutionary} for estimating the evolutionary spectra
\ Here, we consider the estimation of $g_{0}(\omega)$ assuming \ that the
process is spatially and temporally stationary. From the above relation, we
obtain by inverting $g_{0}(\omega)=\frac{1}{2\pi}\sum e^{-iu\omega}c(0,u)$. We
further note$\ $that if the process is fully symmetric (see
\cite{gneiting2002nonseparable}) then $c\left(  \mathbf{h},u\right)  =c\left(
-\mathbf{h},-u\right)  $ and $f\left(  \underline{\lambda},\omega\right)
=f\left(  -\underline{\lambda},-\omega\right)  $ and $f\left(  \underline
{\lambda},\omega\right)  >0$ for all $\underline{\lambda}$ and $\omega$. Here
$\underline{\lambda}$ is the frequency associated with spatial coordinates
(usually called the wave number) and $\omega$ is the temporal frequency. In
the following we define the discrete Fourier Transforms of the stationary
process and summarise in the appendix their well known properties which will
be used (for details refer to \ \cite{Brill-book-01}).

Let $\left\{  Y_{t}\left(  \mathbf{s}_{i}\right)  |\,i=1,2,\ldots
m;\ t=1,\ldots,n\right\}  $ be a sample from the zero mean spatio-temporal
stationary process $\left\{  Y_{t}\left(  \mathbf{s}\right)  \right\}
$.\ Consider the time series data \\ $\left\{  Y_{t}\left(  \mathbf{s}%
_{i}\right)  |\,\ t=1,\ldots,n\right\}  $ at the location $\mathbf{s}_{i}$,
and define the Discrete Fourier transform (DFT)
\[
J_{\mathbf{s}_{i}}\left(  \omega_{k}\right)  =\frac{1}{\sqrt{2\pi n}}%
\sum\limits_{t=1}^{n}Y_{t}\left(  \mathbf{s}_{i}\right)  e^{-it\omega_{k}},
\]
where $\omega_{k}=\frac{2\pi k}{n},\ \ k=0,\ \pm1,\ \ldots,\ \pm\left[
\frac{n}{2}\right]  $. In practice one uses Fast Fourier Transform algorithm
to compute the DFT. From the above, by inversion we get%
\[
Y_{t}\left(  \mathbf{s}\right)  =\sqrt{\frac{n}{2\pi}}\int\limits_{-\pi}^{\pi
}J_{\mathbf{s}}\left(  \omega\right)  e^{it\omega}d\omega\text{.}%
\]
The above integral representation shows that the process can be decomposed
into various sine and cosine terms and complex valued DFT's as the amplitudes.
We also see from the above that there is a one to one correspondence between
the DFT's and the data, a property we use later for prediction.

\section{A \ Complex Stochastic Partial Differential Equation\ (CSPDE) and an
Expression for the spectrum $g_{||\mathbf{h}||}(\omega)$}

In the following we consider the spatial \ model proposed by
\cite{Whittle1954}, spatio-temporal models proposed by \cite{jones1997models},
\cite{lindgren2011explicit}, and \cite{sigrist2015stochastic} and derive the
spectral properties of the processes satisfying these models. Our object here
is to define an alternative frequency domain based complex stochastic partial
differential equation for the discrete Fourier transforms \ and \ derive
expressions for the spectrum and the covariance function of the process
satisfying the CSPDE model. We show that the spectra obtained from the
stochastic partial differential equations defined by \cite{Whittle1954},
\cite{jones1997models}, \cite{lindgren2011explicit},
\cite{sigrist2015stochastic} can be derived as special cases of the CSPDE
model defined here.

It is well known that to study turbulence, dissipation of heat or fluid,
equations such as Laplace equation, Heat equation, Wave equation and Helmholtz
equation are often used. Laplace equation is used to describe the static
behavior of the material (say fluid) where as Heat and Wave equations are used
to describe the dynamic behavior and are usually called Diffusion equations.
Stochastic version of the Laplace equation was used by \cite{Whittle1954} to
study the correlation pattern of soil fertility in agricultural \ uniformity
trials at different locations \ and by studying the solution of the Stochastic
Laplace equation, \cite{Whittle1954} has shown that the correlation of the
yields at points at '$s$' units apart falls off\ \ as a power of '$s^{-1}$', a
property observed by agricultural scientists.\ We briefly discuss the models
proposed by the above researchers.

\subsection{Stochastic version of the Laplace Equation\newline%
(\cite{Whittle1954})}

\ For illustration purposes, let us assume $d=2.$ Let $Y\left(  \mathbf{s}%
\right)  $ (here the spatial coordinates are denoted by $\mathbf{s}%
=(s_{1},\ s_{2})$) denote the stationary spatial random field. In the case of
the example considered by \cite{Whittle1954}, $Y(\mathbf{s})$ denotes the
yield at the location $s.$ Let $\ \bigtriangledown=\frac{\partial^{2}%
}{\partial s_{1}^{2}}+\frac{\partial^{2}}{\partial s_{2}^{2}}$ be the Laplace
operator. \cite{Whittle1954} defined the model $(\bigtriangledown-\gamma
^{2})Y\left(  \mathbf{s}\right)  =e(\mathbf{s})$, where $e(\mathbf{s})$ is
defined \ as spatial \ Gaussian white noise \ and $\gamma$ is a scale
parameter. We can obtain an expression for the spectrum of the process
satisfying the above model by considering the spectral representation of the
process given by\ $Y(\mathbf{s})=\int e^{i\mathbf{s}.\underline{\lambda}}%
$\ $dZ_{y}(\underline{\lambda})$, where $Z_{Y}(\underline{\lambda})$ is an
orthogonal function with \ $E(dZ_{Y}(\underline{\lambda}))=0,\ $and
$E|dZ_{Y}(\underline{\lambda})|^{2}=f_{Y}(\underline{\lambda})d\underline
{\lambda}$. Here $\underline{\lambda}=(\lambda_{1},\ \lambda_{2})$ corresponds
to the spatial frequency (known as wave number), $f_{Y}(\underline{\lambda})$
is defined as the spatial spectrum. \ We can define a similar spectral
representation for Gaussian white noise process $e(\mathbf{s})$ with
$Z_{e}(\underline{\lambda})$ denoting the orthogonal random set function \ of
the process. \ By substituting the spectral representations for $Y(\mathbf{s}%
)$ and $e(\mathbf{s})$ and equating the integrands and taking expectations of
the modulus squares both sides, \ we can show that the spectral density
function of $Y(\mathbf{s})$ is given by $f_{Y}(\underline{\lambda}%
)\propto1/(\lambda_{1}^{2}+\lambda_{2}^{2}+\gamma^{2})^{2}$.

If the process is isotropic (see \cite{SteinWeiss}, Ch. IV. Theorem 1.1.), we
can show by inversion,that the corresponding \ spatial covariance function at
Euclidean distance $||\mathbf{h}||$ is given by $(||\mathbf{h}||/\gamma
)K_{1}(\gamma||\mathbf{h}||)$,$\ $ where $K_{1}(.)$ is the modified Bessel
function of the second kind of the first order. The covariance function
obtained belongs to Matern Class of covariance functions. The above model is a
static version of the dynamic model \ considered below.

\subsection{Stochastic version of the Heat Equation (\cite{jones1997models})}

Now consider the model $[\frac{\partial}{\partial t}+(\bigtriangledown
-\gamma^{2})]$ $Y_{t}(\mathbf{s})=e_{t}(\mathbf{s),\ }$where the white noise
process $\{e_{t}(\mathbf{s)\}}$ is defined \ as above.\texttt{\ }The models
considered by \cite{lindgren2011explicit} and recently by
\cite{sigrist2015stochastic} are variations of the above diffusion model. If
we set transport direction vector (see \cite{sigrist2015stochastic}) zero and
diffusion matrix to identity matrix in the models by
\cite{lindgren2011explicit} and \cite{sigrist2015stochastic}, we get the model
defined by \cite{jones1997models}. By substituting the spectral
representations\texttt{\ }for the processes\texttt{\ }$\{Y_{t}(\mathbf{s})\}$
and $\{e_{t}(\mathbf{s)\}}$ and equating the integrands, and after taking
expectations we can show that the spatio-temporal spectral density function of
the process $\{Y_{t}(\mathbf{s})\}$ is given by\ \ $f_{Y}\left(
\underline{\lambda},\omega\right)  $\ \ $\propto1/[((\lambda_{1}^{2}%
+\lambda_{2}^{2}+\gamma^{2})^{2}+\omega^{2}]$.\

\subsection{Stochastic version of the Wave Equation}

This is a dynamic stochastic version of the classical Wave equation used to
describe sound waves, water waves, light waves arising in fields like
acoustics, fluid dynamics etc.. We are interested in the statistical
properties of the process. Consider the model\ $[\frac{\partial^{2}}{\partial
t^{2}}+(\bigtriangledown-\gamma^{2})]$ $Y_{t}(s)=e_{t}(\mathbf{s)} $\textbf{.
}By substituting the spectral representations, and taking expectations, we can
show that the spatio-temporal spectrum is given by\ \ $f_{Y}\left(
\underline{\lambda},\omega\right)  \propto1/[(\lambda_{1}^{2}+\lambda_{2}%
^{2}+\gamma^{2})+\omega^{2}]^{2}$.\

\subsection{CSPDE and an expression for the spectrum $g_{||\mathbf{h}%
||}(\omega)$}

We note that the above equations considered by \cite{jones1997models},
\cite{lindgren2011explicit}, and \cite{sigrist2015stochastic}, include a first
order time derivative $\frac{\partial}{\partial t}$ in the operators. This is
equivalent to assuming that the temporal dynamics \ in the spatio-temporal
process can be explained by\ an autoregressive model of order one and in a
similar way the inclusion of the second order time derivative in the wave
equation is equivalent to assuming that the temporal dynamics can be explained
by an autoregressive model of order two. These specific assumptions can be
unrealistic in some situations. \ In view of this, we propose \ a model which
includes a nonparametric function which is a polynomial in $e^{-i\omega}$, and
by including this function in our operator, we can derive the spectra defined
by the processes satisfying the above models as special cases. To arrive at
the model, let us consider once again the Laplace\ operator
\ $(\bigtriangledown-\gamma^{2})$ operating on the process $\{Y_{t}%
(\mathbf{s})\}$. We have shown (see section 2)
\[
Y_{t}\left(  \mathbf{s}\right)  =\sqrt{\frac{n}{2\pi}}\int\limits_{-\pi}^{\pi
}J_{\mathbf{s}}\left(  \omega\right)  e^{it\omega}d\omega\text{.}%
\]
Multiplying$\ $both sides of the above equation by the operator
$(\bigtriangledown-\gamma^{2})$, we get \ $(\bigtriangledown-\gamma^{2}%
)Y_{t}\left(  \mathbf{s}\right)  =\sqrt{\frac{n}{2\pi}}\int\limits_{-\pi}%
^{\pi}(\bigtriangledown-\gamma^{2})J_{\mathbf{s}}\left(  \omega\right)
e^{it\omega}d\omega$. This relation shows that operating on the process
$Y_{t}(\mathbf{s})$ is equivalent to operating on the complex valued DFT
$J_{\mathbf{s}}\left(  \omega\right)  $ at a fixed frequency $\omega$ and then
integrating over all the frequencies. In other words, just like\ the
interpretation we have for the spectral representation which \ is a frequency
decomposition of the process in terms of sine and cosine functions and the
contribution of each frequency is measured by the corresponding amplitude
$J_{\mathbf{s}}\left(  \omega\right)  $, the above frequency domain
\ modelling is equivalent to \ modelling the complex valued process (DFT) for
each frequency $\omega$.\ Later we will obtain an expression for the
covariance function of the complex valued process $\{J_{\mathbf{s}}\left(
\omega\right)  \}$ satisfying CSPDE which will be in terms of the temporal
spectrum $g_{0}(\omega)$ and the spatial distance $||h||$. We will show this
covariance function is necessary for spatio-temporal prediction. In the course
of the derivation of this result, we will also obtain an expression for the
second order spectrum \ which is a function of the spatial \ frequency (wave
number) and the temporal frequency. The obtained spectrum is strictly greater
than zero implying that the corresponding spatio-temporal covariance function
is positive definite.\ Also the spectrum \ obtained is non separable.

We now show \ that the above models can be considered as special cases of the
following CSPDE model.

\begin{lemma}
Let d=2. Consider the complex stochastic partial differential equation
\[
\left[  \frac{\partial^{2}}{\partial s_{1}^{2}}+\frac{\partial^{2}}{\partial
s_{2}^{2}}-\gamma\left(  \omega\right)  \right]  J_{\mathbf{s}}\left(
\omega\right)  =J_{\mathbf{s},e}\left(  \omega\right)  ,
\]
where $\gamma\left(  \omega\right)  $\ is a complex valued function. Let
\ $\gamma\left(  \omega\right)  $\ $=c(\omega)+i$ $b(\omega).$ Then the
spectral density function of the process $\{Y_{t}\left(  \mathbf{s}\right)
\}$\ satisfying the above model is given by\
\[
f_{Y}\left(  \underline{\lambda},\omega\right)  \ \propto\frac{1}%
{[((\lambda_{1}^{2}+\lambda_{2}^{2}+c(\omega))^{2}+b^{2}(\omega)]}%
\]
\ \
\end{lemma}

\begin{proof}
We can proceed as before to obtain the above expression and hence the\ proof
is omitted.
\end{proof}

Let us now consider the following special cases.\

\begin{enumerate}
\item Let $c(\omega)=\gamma^{2}$, $b(\omega)=\omega$. Substitute these in the
equation $\ f_{Y}\left(  \underline{\lambda},\omega\right)  $ given above, we
obtain the spectrum of \ the process satisfying the \ Heat\ equation
considered by \cite{jones1997models}, \cite{lindgren2011explicit} and
\cite{sigrist2015stochastic}. \

\item To obtain the spectrum of the Wave equation, let $\gamma\left(
\omega\right)  $ be real valued, and let $\gamma\left(  \omega\right)
=c(\omega)=\gamma^{2}+\omega^{2}.$ Substitution of this gives us the spectrum
corresponding to the Wave equation.
\end{enumerate}

Through the above examples we have shown that by defining the \ stochastic
version of the Laplacian \ model in terms of the frequency dependent
\ nonparametric function $\gamma\left(  \omega\right)  $, the second order
properties of the classical equations can be obtained as special cases.

We will now state the \ main model \ and derive expressions for the spectrum
and for the covariance function$\ g_{\parallel\mathbf{h}\parallel}(\omega) $
which are functions of the Euclidean distance $\left\Vert \mathbf{h}%
\right\Vert $ \ and temporal spectral frequency $\omega$. We will state the
results for $d=2$\ and later consider its generalization for all $d$.

\begin{theorem}
Let $J_{\mathbf{s}}\left(  \omega\right)  $ be the discrete Fourier transform
of the data $\{Y_{t}(\mathbf{s})|t=1,\ 2,\ ...n\}$ at the location
$\mathbf{s}.$ Let $\upsilon$ $>$\ $0$, and let $\{J_{\mathbf{s}}(\omega)\}$
satisfy the model%
\begin{equation}
\left[  \frac{\partial^{2}}{\partial s_{1}^{2}}+\frac{\partial^{2}}{\partial
s_{2}^{2}}-\left\vert c\left(  \omega\right)  \right\vert ^{2}\right]  ^{\nu
}J_{\mathbf{s}}\left(  \omega\right)  =J_{\mathbf{s},e}\left(  \omega\right)
, \label{eq 3.4.}%
\end{equation}
where $J_{\mathbf{s}}\left(  \omega\right)  $ and $J_{\mathbf{s},e}\left(
\omega\right)  $ are given by (\ref{eq 2.8.}) and (\ref{DFT_WN}). Then the
second order spectral density function $f_{Y}\left(  \underline{\lambda
},\omega\right)  $ of the process$\{Y_{t}(\mathbf{s})\}$ is given by%
\[
f_{Y}\left(  \underline{\lambda},\omega\right)  =\frac{\sigma_{e}^{2}}{\left(
2\pi\right)  ^{2}\left(  \lambda_{1}^{2}+\lambda_{2}^{2}+\left\vert c\left(
\omega\right)  \right\vert ^{2}\right)  ^{2\nu}}.
\]
If the stationary spatio-temporal process is isotropic, then the covariance
function between the discrete Fourier Transforms $J_{\mathbf{s}}\left(
\omega\right)  $ and $J_{\mathbf{s+h}}\left(  \omega\right)  $ is given by%
\[
g_{\parallel\mathbf{h}\parallel}(\omega)=Cov(J_{\mathbf{s}}\left(
\omega\right)  ,J_{\mathbf{s+h}}\left(  \omega\right)  )=\ \frac{\sigma
_{e}^{2}}{2\pi}\left(  \frac{\left\Vert \mathbf{h}\right\Vert }{2\left\vert
c\left(  \omega\right)  \right\vert }\right)  ^{2\nu-1}\frac{K_{2\nu-1}\left(
\left\vert c\left(  \omega\right)  \right\vert \left\Vert \mathbf{h}%
\right\Vert \right)  }{\Gamma(2\nu)},
\]
where $\left\Vert \mathbf{h}\right\Vert =\left(  h_{1}^{2}+h_{2}^{2}\right)
^{\frac{1}{2}}$ and $K_{\nu}\left(  x\right)  $ is the modified Bessel
function of the second kind of order $\nu$.
\end{theorem}

Later we will see the significance of the frequency dependent function
$c\left(  \omega\right)  $ in the above model.

\begin{proof}
Substitute the spectral representations for $J_{\mathbf{s}}\left(
\omega\right)  $ and $J_{\mathbf{s},e}\left(  \omega\right)  $ given by
(\ref{eq 2.8.}) and (\ref{DFT_WN}) \ and taking the operators inside the
integrands and equating the integrands both sides of the equation (this is
valid because of the uniqueness of the Fourier transforms), we obtain
\begin{equation}
\left(  -\lambda_{1}^{2}-\lambda_{2}^{2}-\left\vert c\left(  \omega\right)
\right\vert ^{2}\right)  ^{\nu}dZ_{Y}\left(  \underline{\lambda}%
,\omega\right)  =dZ_{e}\left(  \underline{\lambda},\omega\right)  \text{, }
\label{eq 3.5.}%
\end{equation}
where $\underline{\lambda}=\left(  \lambda_{1},\lambda_{2}\right)  $. Taking
the modulus squares, and taking expectations both sides of the\ modulus
squares we obtain the spatio-temporal spectral density function of the
spatio-temporal process $Y_{t}\left(  \mathbf{s}\right)  $ satisfying the
above model (\ref{eq 3.4.}) and it is given by
\begin{equation}
f_{Y}\left(  \underline{\lambda},\omega\right)  =\frac{\sigma_{e}^{2}}{\left(
2\pi\right)  ^{2}\left(  \lambda_{1}^{2}+\lambda_{2}^{2}+\left\vert c\left(
\omega\right)  \right\vert ^{2}\right)  ^{2\nu}}. \label{eq 3.6.}%
\end{equation}
which is the stated result. \newline We note that the above spectral density
is real and strictly positive, and this implies that the associated
\ spatio-temporal covariance function is positive definite. Further the
spectral density function given above belongs to a nonseparable class of
process. Since the spectrum depends on the distance of the wave numbers
$\lambda_{1}^{2}+\lambda_{2}^{2}$ its Fourier transform will depend on the
distance between locations as well (see \cite{SteinWeiss} Ch. IV. Theorem
1.1.), hence the process is isotropic. Now, to obtain the covariance function
$g_{\parallel\mathbf{h}\parallel}(\omega)$, we need to take its inverse
Fourier transform. We use the result used by \cite{Whittle1954} (equation (65)
of the paper)
\[
\frac{1}{4\pi^{2}}\int\int\frac{e^{i\left(  x\omega_{1}+y\omega_{2}\right)  }%
}{\left(  \omega_{1}^{2}+\omega_{2}^{2}+\alpha^{2}\right)  ^{\mu+1}}%
d\omega_{1}d\omega_{2}=\frac{1}{2\pi}\left(  \frac{r}{2\alpha}\right)  ^{\mu
}\frac{K_{\mu}\left(  \alpha r\right)  }{\Gamma\left(  \mu+1\right)  }\text{,
}%
\]
where $r=\left(  x^{2}+y^{2}\right)  ^{\frac{1}{2}}$, $K_{\mu}\left(
x\right)  $ is the modified Bessel function of the second kind of order $\mu
$.\ We use the above result to obtain the inverse transform of $f_{Y}\left(
\underline{\lambda},\omega\right)  $, given by (\ref{eq 3.6.})$.$ Taking the
inverse transform over the wave number $\underline{\lambda}$ only \ (for fixed
temporal frequency $\omega$), we obtain
\begin{align}
g_{\parallel\mathbf{h}\parallel}(\omega)  &  =\frac{\sigma_{e}^{2}}{\left(
2\pi\right)  ^{2}}\int\int\frac{e^{i\left(  h_{1}\lambda_{1}+h_{2}\lambda
_{2}\right)  }}{\left(  \lambda_{1}^{2}+\lambda_{2}^{2}+\left\vert c\left(
\omega\right)  \right\vert ^{2}\right)  ^{2\nu}}d\lambda_{1}d\lambda
_{2}\nonumber\\
&  =\frac{\sigma_{e}^{2}}{2\pi}\left(  \frac{\left\Vert \mathbf{h}\right\Vert
}{2\left\vert c\left(  \omega\right)  \right\vert }\right)  ^{2\nu-1}%
\frac{K_{2\nu-1}\left(  \left\vert c\left(  \omega\right)  \right\vert
\left\Vert \mathbf{h}\right\Vert \right)  }{\Gamma(2\nu)}. \label{eq 3.7.}%
\end{align}

\end{proof}

The above interesting expression shows that the covariance function between
two discrete Fourier Transforms separated by the spatial distance $\left\Vert
\mathbf{h}\right\Vert $ again can be written in terms Matern and Whittle class
of covariance functions.

\textit{However, the most important and fundamental difference between this
expression and other covariance expressions given by other authors is that the
argument of the Bessel function derived above is not only a function of the
spatial distance, but also a function of the frequency dependent scaling
function }$\left\vert c\left(  \omega\right)  \right\vert $\textit{\ \ which
is related to the second order temporal spectral function. This will be shown
in the following lemma.}

To see the significance of inclusion of $|c(\omega)|$ in the model
(\ref{eq 3.4.}), we consider the limiting behavior of $g_{\parallel
\mathbf{h}\parallel}(\omega)$ as $||\mathbf{h}||\rightarrow0$. We have noted
earlier that $Var\left(  J_{\mathbf{s}}\left(  \omega\right)  \right)  $ given
by (\ref{eq 2.2.}) (see Appendix A), is proportional to the spectral density
function $g_{0}\left(  \omega\right)  \ $of the random process for all $s$. So
it is interesting to examine the behavior of $g_{\parallel\mathbf{h}\parallel
}(\omega)$ when $\left\Vert \mathbf{h}\right\Vert \rightarrow0$, as the limit
must tend to the second order spectral density function $g_{0}\left(
\omega\right)  $ of the process $\{Y_{t}\left(  \mathbf{s}\right)  \}$. We
state the result in the following Lemma.

\begin{lemma}
For the above isotropic process, and under the conditions stated above, as
$\left\Vert \mathbf{h}\right\Vert \rightarrow0$, $g_{\parallel\mathbf{h}%
\parallel}(\omega)$ tends to $\ $%
\begin{equation}
g_{0}\left(  \omega\right)  =\frac{\sigma_{e}^{2}}{2\left(  \left\vert
c\left(  \omega\right)  \right\vert ^{2}\right)  ^{2\nu-1}\left(
2\nu-1\right)  }. \label{eq 3.10.}%
\end{equation}

\end{lemma}

\begin{proof}
It is well known that, for all $\nu>0$,
\begin{equation}
\lim_{x\rightarrow0}\frac{x^{\nu}K_{\nu}\left(  x\right)  }{2^{\nu-1}%
\Gamma\left(  \nu\right)  }=1\text{.} \label{eq 3.9.}%
\end{equation}
Therefore, if we take the limit of $g_{\parallel\mathbf{h}\parallel}(\omega)$
given by (\ref{eq 3.7.}) as $\left\Vert \mathbf{h}\right\Vert \rightarrow0$,
we get \ the stated result,
\end{proof}

\subsubsection{Special Case:}

Let us consider the case $\nu=1$. Then\ from the equation (\ref{eq 3.7.}) we
have%
\begin{equation}
g_{\parallel\mathbf{h}\parallel}(\omega)=\frac{\sigma_{e}^{2}}{2\pi}\left(
\frac{\left\Vert \mathbf{h}\right\Vert }{2\left\vert c\left(  \omega\right)
\right\vert }\right)  K_{1}\left(  \left\vert c\left(  \omega\right)
\right\vert \left\Vert \mathbf{h}\right\Vert \right)  , \label{eq 3.7.a}%
\end{equation}
and \ from the equation (\ref{eq 3.10.}) we \ have\
\begin{equation}
g_{0}\left(  \omega\right)  \ =\frac{\sigma_{e}^{2}}{2\left\vert c\left(
\omega\right)  \right\vert ^{2}}, \label{eq 3.7.0}%
\end{equation}
which implies that $\left\vert c\left(  \omega\right)  \right\vert ^{2}$ is
proportional to $g_{0}^{-1}\left(  \omega\right)  $, which is defined as the
inverse second order spectral density function of the process. Let us assume
that $g_{0}^{-1}\left(  \omega\right)  $ is absolutely integrable, then
$g_{0}^{-1}\left(  \omega\right)  $ can be expanded in Fourier series
\[
g_{0}^{-1}\left(  \omega\right)  =\frac{1}{2\pi}\sum\limits_{k=-\infty
}^{\infty}ci\left(  k\right)  \cos k\omega,\left\vert \omega\right\vert
\leq\pi\text{, }%
\]
where we used the fact that $g_{0}^{-1}\left(  \omega\right)  =g_{0}%
^{-1}\left(  -\omega\right)  $. The coefficients $\left\{  ci\left(  k\right)
\right\}  $ are usually known as inverse autocovariances, and sometimes are
used to estimate the orders of the linear time series models. For example, if
the series $\left\{  Y_{t}\left(  \mathbf{s}\right)  \right\}  $ satisfy (for
a given $s$) an autoregressive model of order $p$, say, then it can easily be
shown that $ci\left(  k\right)  =0$ for all $k>p$. In view of this interesting
property one can use the inverse auto-covariances to determine the order of
the \ linear AR models.\ We note further that the covariance function
$g_{\parallel\mathbf{h}\parallel}(\omega)$ given above \ is in terms of the
modified Bessel function, the argument of the Bessel function is a product of
the spatial distance $\left\Vert \mathbf{h}\right\Vert $ \ and the inverse
temporal spectrum $g_{0}^{-1}\left(  \omega\right)  $. Therefore the rate of
convergence of the covariance function to tend to zero as $\left\Vert
\mathbf{h}\right\Vert \rightarrow\infty\ $\ depends on the second order
temporal spectrum of the process at the frequency $\omega$.

From (\ref{eq 3.7.a}) and (\ref{eq 3.7.0}) we can also obtain an expression
for the auto-correlation function. We have the auto-correlation function when
$d=2$, and for all $\upsilon>0,$
\begin{align}
\rho\left(  \left\Vert \mathbf{h}\right\Vert ,\omega\right)   &
=\frac{g_{\parallel\mathbf{h}\parallel}(\omega)}{\ g_{0}\left(  \omega\right)
}\nonumber\\
&  =\frac{\left(  \left\Vert \mathbf{h}\right\Vert \left\vert c\left(
\omega\right)  \right\vert \right)  ^{2\nu-1}}{2^{2\nu-2}\Gamma\left(
2\nu-1\right)  }K_{2\nu-1}\left(  \left\vert c\left(  \omega\right)
\right\vert \left\Vert \mathbf{h}\right\Vert \right)  \text{.}
\label{eq 3.11.}%
\end{align}
It is interesting to note that $\rho\left(  \left\Vert \mathbf{h}\right\Vert
,\omega\right)  $ is in fact the coherency coefficient between two Discrete
Fourier Transforms separated by the spatial distance $\left\Vert
\mathbf{h}\right\Vert $ at \ the frequency $\omega$. We now consider the
generalization of Theorem 1.\ \

\begin{theorem}
Let $\nu>0$ and $\ d\geq2$. Let the Discrete Fourier Transform
$\ J_{\mathbf{s}}\left(  \omega\right)  $ satisfy the equation
\[
\left[  \sum_{i=1}^{d}\frac{\partial^{2}}{\partial s_{i}^{2}}-\left\vert
c\left(  \omega\right)  \right\vert ^{2}\right]  ^{\nu}J_{\mathbf{s}}\left(
\omega\right)  =J_{\mathbf{s},e}\left(  \omega\right)  ,
\]
Then the second order spectral density function is given by$\ $%
\[
f_{Y}\left(  \underline{\lambda},\omega\right)  =\frac{\sigma_{e}^{2}}{\left(
2\pi\right)  ^{d}}\frac{1}{\left(  \sum\lambda_{i}^{2}+\left\vert c\left(
\omega\right)  \right\vert ^{2}\right)  ^{2\nu}}.
\]
If the process is isotropic then the covariance function is given by\
\begin{align*}
g_{\parallel\mathbf{h}\parallel}(\omega)\  &  =\ Cov(\ J_{\mathbf{s}}\left(
\omega\right)  ,J_{\mathbf{s+h}}\left(  \omega\right)  )\ \\
&  =\frac{\sigma_{e}^{2}}{\left(  2\pi\right)  ^{\frac{d}{2}}2^{2\nu-1}%
\Gamma\left(  2\nu\right)  }\left(  \frac{\left\Vert \mathbf{h}\right\Vert
}{\left\vert c\left(  \omega\right)  \right\vert }\right)  ^{2\nu-\frac{d}{2}%
}K_{2\nu-\frac{d}{2}}\left(  \left\Vert \mathbf{h}\right\Vert \left\vert
c\left(  \omega\right)  \right\vert \right)  \ .
\end{align*}

\end{theorem}

\begin{proof}
By proceeding as in Theorem 1, we can show that the spectral density function
is given by
\[
\ f_{Y}\left(  \underline{\lambda},\omega\right)  =\frac{\sigma_{e}^{2}%
}{\left(  2\pi\right)  ^{d}}\frac{1}{\left(  \sum\lambda_{i}^{2}+\left\vert
c\left(  \omega\right)  \right\vert ^{2}\right)  ^{2\nu}}.
\]
To obtain the inverse transform we proceed as follows. Let $\rho=\left\Vert
\underline{\lambda}\right\Vert $. We have,
\begin{align*}
g_{\parallel\mathbf{h}\parallel}\left(  \omega\right)  \  &  =\frac{\sigma
_{e}^{2}}{\left(  2\pi\right)  ^{d}}\int_{\mathbb{R}^{d}}\frac
{e^{-i\mathbf{h\cdot}\underline{\lambda}}}{\left(  \left\Vert \underline
{\lambda}\right\Vert ^{2}+\left\vert c\left(  \omega\right)  \right\vert
^{2}\right)  ^{2\nu}}d\underline{\lambda}\\
&  =\frac{\sigma_{e}^{2}}{\left(  2\pi\right)  ^{d}}\int_{0}^{\infty}%
\frac{\rho^{d-1}}{\left(  \rho^{2}+\left\vert c\left(  \omega\right)
\right\vert ^{2}\right)  ^{2\nu}}\int_{\mathbb{S}_{d-1}}e^{-i\rho\left\Vert
\mathbf{h}\right\Vert \cos\alpha}d\Omega d\rho,
\end{align*}
where $\mathbb{S}_{d-1}$ is the unit sphere in $\mathbb{R}^{d}$ and $\Omega$
is Lebesgue element of surface area on $\mathbb{S}_{d-1}$. We know further
\[
\int_{\mathbb{S}_{d-1}}e^{-i\rho\left\Vert \mathbf{h}\right\Vert \cos\alpha
}d\Omega=\left(  2\pi\right)  ^{\frac{d}{2}}\left(  \rho\left\Vert
\mathbf{h}\right\Vert \right)  ^{-\frac{d}{2}+1}\mathcal{J}_{\frac{d}{2}%
-1}\left(  \rho\left\Vert \mathbf{h}\right\Vert \right)  ,
\]
where $\mathcal{J}_{\frac{d}{2}-1}$ denotes the Bessel function of the first
kind, see \cite{SteinWeiss}, p. 176. Now we use Hankel-Nicholson Type
Integral, see \cite{Abramowit12}, 11.4.44, if $d<4\nu+3$, then
\[
\int_{0}^{\infty}\frac{\mathcal{J}_{\frac{d}{2}-1}\left(  r\rho\right)
}{\left(  \rho^{2}+\left\vert c\left(  \omega\right)  \right\vert ^{2}\right)
^{2\nu}}\rho^{\frac{d}{2}}d\rho=\frac{r^{2\nu-1}\left\vert c\left(
\omega\right)  \right\vert ^{\frac{d}{2}-2\nu}}{2^{2\nu-1}\Gamma\left(
2\nu\right)  }K_{\frac{d}{2}-2\nu}\left(  r\left\vert c\left(  \omega\right)
\right\vert \right)  .
\]
Using the above integrals and noting $K_{\frac{d}{2}-2\nu}=K_{2\nu-\frac{d}%
{2}}$, for all $d,\ $ the covariance function can be shown to be
\[
g_{\parallel\mathbf{h}\parallel}(\omega)=\frac{\sigma_{e}^{2}}{\left(
2\pi\right)  ^{\frac{d}{2}}2^{2\nu-1}\Gamma\left(  2\nu\right)  }\left(
\frac{\left\Vert \mathbf{h}\right\Vert }{\left\vert c\left(  \omega\right)
\right\vert }\right)  ^{2\nu-\frac{d}{2}}K_{2\nu-\frac{d}{2}}\left(
\left\Vert \mathbf{h}\right\Vert \left\vert c\left(  \omega\right)
\right\vert \right)  ,
\]
and the auto-correlation function is
\[
\rho\left(  \left\Vert \mathbf{h}\right\Vert ,\omega\right)  =\frac{\left(
\left\Vert \mathbf{h}\right\Vert \left\vert c\left(  \omega\right)
\right\vert \right)  ^{2\nu-\frac{d}{2}}}{2^{2\nu-\frac{d}{2}-1}\Gamma\left(
2\nu-\frac{d}{2}\right)  }K_{2\nu-\frac{d}{2}}\left(  \left\Vert
\mathbf{h}\right\Vert \left\vert c\left(  \omega\right)  \right\vert \right)
,
\]
since
\[
\ g_{0}\left(  \omega\right)  =\frac{\sigma_{e}^{2}}{\left(  2\pi\right)
^{\frac{d}{2}}2^{\frac{d}{2}}\left(  \left\vert c\left(  \omega\right)
\right\vert ^{2}\right)  ^{2\nu-\frac{d}{2}}}\frac{\Gamma\left(  2\nu-\frac
{d}{2}\right)  }{\Gamma\left(  2\nu\right)  }.
\]

\end{proof}

In the following we consider prediction of the data, and the optimal
predictive function \ is given in terms of the Discrete Fourier Transforms and
the covariance function $g_{\parallel\mathbf{h}\parallel}(\omega)$.

\section{Spatio-temporal Prediction}

Our object in this section is to estimate $\left\{  Y_{t}\left(
\mathbf{s}\right)  |\,t=1,2,\ldots,n\right\}  $ at \ the location
$\mathbf{s}_{0}$ given the \ $m$ time series $\left\{  Y_{t}\left(
\mathbf{s}_{i}\right)  |\ i=1,2,\ldots,m;t=1,2,\ldots,n\right\}  $ from the
spatio-temporal stationary, isotropic process $\left\{  Y_{t}\left(
\mathbf{s}\right)  \right\}  $. In other words, we are estimating the entire
data set at \ the location $\mathbf{s}_{0}.$ \ Using the estimated
observations \ at the location $\mathbf{s}_{0}$, we can also obtain the
optimal linear predictors for the future values at the location $\mathbf{s}%
_{0}$. As in the case of the observed data $\left\{  Y_{t}\left(
\mathbf{s}_{i}\right)  \right\}  $, we define the discrete Fourier transform
$J_{\mathbf{s}_{0}}\left(  \omega\right)  $ of $\left\{  Y_{t}\left(
\mathbf{s}_{0}\right)  \right\}  ,\ $ and estimate the Fourier transform
$J_{\mathbf{s}_{0}}\left(  \omega\right)  $ for all $\omega$. Using the
inverse Fourier Transform, we recover the data \ $Y_{t}\left(  \mathbf{s}%
_{0}\right)  $ for all $1\leq t$\ $\leq n.$ We pointed out earlier that there
is a one to one correspondence between the discrete Fourier Transforms and the
data. We have shown\ earlier, if
\begin{equation}
J_{\mathbf{s}_{0}}\left(  \omega\right)  =\frac{1}{\sqrt{(2\pi}n)}%
\sum\limits_{t=1}^{n}Y_{t}\left(  \mathbf{s}_{0}\right)  e^{-it\omega}\text{,
}\label{eq 5.1.}%
\end{equation}
then we have
\begin{equation}
Y_{t}\left(  \mathbf{s}_{0}\right)  =\sqrt{\frac{n}{2\pi}}\int\limits_{-\pi
}^{\pi}J_{\mathbf{s}_{0}}\left(  \omega\right)  e^{it\omega}d\omega
\text{.}\label{eq 5.1. a}%
\end{equation}
Consider the vector of the discrete Fourier transforms obtained from all \ the
$m$ locations at the frequency $\omega,\ $
\[
\underline{J}_{m}^{\prime}\left(  \omega\right)  =\left[  J_{\mathbf{s}_{1}%
}\left(  \omega\right)  ,J_{\mathbf{s}_{2}}\left(  \omega\right)
,\ldots,J_{\mathbf{s}_{m}}\left(  \omega\right)  \right]  \text{.}%
\]
We note that
\begin{align}
E\left[  \underline{J}_{m}\left(  \omega\right)  \right]   &  =0,\nonumber\\
E\left[  \underline{J}_{m}\left(  \omega\right)  \underline{J}_{m}^{\ast
}\left(  \omega\right)  \right]   &  =F_{m}\left(  \omega\right)  \text{,
}\label{eq 5.2.}%
\end{align}
where$\ $the\ real, symmetric, positive definite square matrix \newline%
$F_{m}\left(  \omega\right)  =\left(  g_{\parallel s_{i}-s_{j}\parallel
}(\omega);\ i,j=1,2,\ldots,m\right)  $, and each\ element $g_{\parallel
s_{i}-s_{j}\parallel}(\omega)$ of $\ $the matrix $F_{m}\left(  \omega\right)
$ is given by (\ref{eq 3.7.}).\ The complex random vector $\underline{J}%
_{m}\left(  \omega\right)  $ has a multivariate complex Gaussian distribution
with mean zero and variance covariance matrix \ $F_{m}\left(  \omega\right)
$. Consider now the $\left(  m+1\right)  $ dimensional \ complex valued random
vector,
\[
\underline{J}_{m+1}^{\prime}\left(  \omega\right)  =\left[  J_{\mathbf{s}%
_{_{o}}}\left(  \omega\right)  ,\underline{J}_{m}^{\prime}\left(
\omega\right)  \right]  .\text{ }%
\]
It can be shown that the mean of the vector is zero, and the variance
covariance matrix is given by
\begin{align*}
E\left[  \underline{J}_{m+1}\left(  \omega\right)  \underline{J}_{m+1}^{\ast
}\left(  \omega\right)  \right]   &  =\left[
\begin{array}
[c]{cc}%
E(J_{s_{_{o}}}\left(  \omega\right)  J_{s_{_{o}}}^{\ast}\left(  \omega\right)
) & E\left(  J_{0}\left(  \omega\right)  \underline{J}_{m}^{\ast\prime}\left(
\omega\right)  \right)  \\
E\left(  \underline{J}_{m}\left(  \omega\right)  J_{0}^{\ast}\left(
\omega\right)  \right)   & E\left(  \underline{J}_{m}\left(  \omega\right)
\underline{J}_{m}^{\ast}\left(  \omega\right)  \right)
\end{array}
\right]  \\
&  =\left[
\begin{array}
[c]{cc}%
g_{0}(\omega) & \underline{G}_{0}^{\prime}\left(  \omega\right)  \\
\underline{G}_{0}\left(  \omega\right)   & F_{m}\left(  \omega\right)
\end{array}
\right]  \text{, }%
\end{align*}
where $g_{0}(\omega)$ is the second order spectral density function of the
spatial process $\ \{Y_{t}\left(  \mathbf{s}_{0}\right)  \}$ and the row
vector $\underline{G}_{0}^{\prime}\left(  \omega\right)  $ is given by
\begin{align*}
\underline{G}_{0}^{\prime}\left(  \omega\right)   &  =E\left[  J_{\mathbf{s}%
_{o}}\left(  \omega\right)  J_{m}^{\ast\prime}\left(  \omega\right)  \right]
\\
&  =[g_{\parallel s_{0}-s_{1}\parallel}(\omega),g_{\parallel s_{0}%
-s_{2}\parallel}(\omega),\ldots,g_{\parallel s_{0}-s_{n}\parallel}(\omega)],
\end{align*}
and $F_{m}\left(  \omega\right)  $ is defined above. Therefore, the optimal
linear least squares predictor of $J_{0}\left(  \omega\right)  $ given the
vector $\underline{J}_{m}\left(  \omega\right)  $, is given by the conditional
expectation
\begin{equation}
E\left[  J_{\mathbf{s}_{o}}\left(  \omega\right)  |\,\underline{J}_{m}\left(
\omega\right)  \right]  =\underline{G}_{0}^{\prime}\left(  \omega\right)
F_{m}^{-1}\left(  \omega\right)  \underline{J}_{m}\left(  \omega\right)
,\label{eq 5.3.}%
\end{equation}
and the minimum mean square prediction error is given by
\begin{equation}
\sigma_{m}^{2}\left(  \omega\right)  =g_{0}(\omega)-\underline{G}_{0}^{\prime
}\left(  \omega\right)  F_{m}^{-1}\left(  \omega\right)  \underline{G}%
_{0}\left(  \omega\right)  \text{.}\label{eq 5.4.}%
\end{equation}

To estimate the data $Y_{t}\left(  \mathbf{s}_{0}\right)  $ for all $t$, we
use the inverse transform (\ref{eq 5.1. a}) and as an estimate of
$J_{\mathbf{s}_{0}}(\omega)$ the right hand expression of (\ref{eq 5.3.}),
where we replace the elements of the matrices $\underline{G}_{0}\left(
\omega\right)  $ and $F_{m}\left(  \omega\right)  $ by their estimates and
obtain
\[
\widehat{Y}_{t}\left(  \mathbf{s}_{0}\right)  =\sqrt{\frac{n}{2\pi}}%
\int\limits_{-\pi}^{\pi}e^{it\omega}\underline{\widehat{G}}_{0}^{\prime
}\left(  \omega\right)  \widehat{F}_{m}^{-1}\left(  \omega\right)
\underline{J}_{m}\left(  \omega\right)  d\omega.
\]

We note $E(Y_{t}(\mathbf{s}_{0}))=0$ \ and $Var(Y_{t}\left(  \mathbf{s}%
_{0}\right)  )\simeq\int\limits_{-\pi}^{\pi}G_{0}^{^{\prime}}\left(
\omega\right)  F_{m}^{-1}\left(  \omega\right)  G_{0}(\omega)d\omega$. We can
show by an application of Parseval's Theorem
\begin{align}
E\left(  Y_{t}\left(  \mathbf{s}_{0}\right)  -\widehat{Y}_{t}\left(
\mathbf{s}_{0}\right)  \right)  ^{2}  &  =E\left\vert \mathcal{F}^{-1}\left(
J_{\mathbf{s}_{o}}\left(  \omega\right)  -\widehat{J}_{\mathbf{s}_{o}}\left(
\omega\right)  \right)  \right\vert ^{2}\nonumber\\
&  =\int\limits_{-\pi}^{\pi}\sigma_{m}^{2}\left(  \omega\right)  d\omega.
\label{StandErr}%
\end{align}
In practice, the above integrals are approximated by finite sums of the form%
\[
\widehat{Y}_{t}\left(  \mathbf{s}_{0}\right)  =\sqrt{\frac{2\pi}{n}}\sum
e^{it\omega_{j}}\underline{\widehat{G}}_{0}^{\prime}\left(  \omega_{j}\right)
\widehat{F}_{m}^{-1}\left(  \omega_{j}\right)  \underline{J}_{m}\left(
\omega_{j}\right)  .
\]
for all $t=1,\ 2,\ ...n,\ $where the estimates $\underline{\widehat{G}}%
_{0}\left(  \omega_{j}\right)  $ and $\widehat{F}_{m}\left(  \omega\right)  $
are substituted for $\underline{G}_{0}\left(  \omega\right)  $ and
$F_{m}\left(  \omega\right)  $ respectively. \ As noted earlier, the vector
$\underline{G}_{0}\left(  \omega\right)  $ and the matrix $F_{m}\left(
\omega\right)  $ have covariance functions $g_{\parallel s_{i}-s_{j}\parallel
}(\omega)$ as their elements. The covariance functions are functions of some
unknown parameters which are related to the spatial correlation and temporal
correlation. From the expression of the covariance function (\ref{eq 3.7.}),
we \ see \ that the parameters to be estimated are $\sigma_{e}^{2}$ (the
variance of the white noise process $e(\mathbf{s},t)$), and the parameters of
the spatio-temporal spectrum $g_{0}(\omega)$. Let us denote the parameter
vector which characterizes $g_{0}(\omega)$ by $\vartheta_{1}$ and the entire
parameter vector by $\vartheta=(\sigma_{e}^{2},\vartheta_{1})$. The parameter
$\nu$ is related to \ the smoothness of the process. In practice one considers
several possible choices for $\nu(>0)$ a priori. The widely used choice is
$\nu=1$. The estimation of the parameter vector $\vartheta$ of the
spatio-temporal covariance function $g_{\parallel\mathbf{h}\parallel}(\omega)$
is extremely important \ and this will be considered in the following section.
We \ now make some comments on computational aspects.

It is interesting and important to note that from the equations (\ref{eq 5.3.}%
) and \ (\ref{eq 5.4.}) that the evaluation of the conditional \ expectation
and the calculation of the minimum mean square error requires inversion of
$\ m\times m$ dimensional matrices (where $m$ corresponds to the number of
locations)only, unlike in the case of time domain approach for prediction
where one needs to invert $\ mn\times mn$ dimensional matrices. In many real
data analysis usually the number of time points $n$ will be \ very large (and
$m$ can be large too). Besides, there is no ordering problem involved here
(see \cite{Cressie2011} p. 324). Once \ we have an expression for \ the
covariance function $g_{\parallel\mathbf{h}\parallel}(\omega)$, all the
elements of the column vector $G_{0}\left(  \omega\right)  $ and the elements
of $F_{m}\left(  \omega\right)  $ are known. By substituting the \ relevant
expressions (or their estimates), we can evaluate (\ref{eq 5.3.}) and
(\ref{eq 5.4.}).

It may be pointed out that there are other approaches for obtaining forecasts
in the context of spatio-temporal data. \cite{sahu2005bayesian} used Bayesian
approach based on MCMC, \cite{ruiz2012new} and \cite{giraldo2010continuous}
\ based their methodology on the assumption that the spatio-temporal data is
of functional data type. \cite{giraldo2010continuous} assumed that the process
can be expanded in terms of some chosen deterministic basis functions with
random coefficients \ and the predictor can also be written as a linear
combination of the same basis functions and the same number of terms. The
solution depends on inversion of \ matrices whose dimensions depend not only
on number of locations and also on the number of Basis functions included in
the expansion of the process and the estimator proposed. All the above
approaches are time domain approaches, and we refer to their papers and papers
there in for more details.

As pointed out earlier, the computation of the predictor depends on the
knowledge of $\underline{G}_{0}\left(  \omega\right)  $ and $F_{m}\left(
\omega\right)  $ which in turn depends on several parameters of the covariance
function $g_{\parallel\mathbf{h}\parallel}(\omega)$. In the following section
we will consider the estimation of the parameters. \ The estimation \ is based
\ on Frequency Variogram approach recently proposed by \cite{SubbaRaoa}. In
their paper, \cite{SubbaRaoa} \ considered the estimation of the parameters of
the covariance function, their asymptotic properties and their efficiency
compared to Gaussian likelihood approach. To avoid repetition, we refer to
\ \cite{SubbaRaoa} for full details.

\section{Estimation of the Parameters of the Covariance function\newline%
$g_{\parallel\mathbf{h}\parallel}(\omega)$ by Frequency Variogram (FV) Method}

We now consider the estimation of the parameters of the covariance function
$g_{\parallel\mathbf{h}\parallel}(\omega)$ using the Frequency variogram
approach recently suggested by \cite{SubbaRaoa}. Here we discuss briefly the
\ FV methodology, and for details, we refer to \cite{SubbaRaoa}. Let
\ $g_{\parallel\mathbf{h}\parallel}(\omega)$\ $=$\ $Cov($\ $J_{\mathbf{s}%
}\left(  \omega\right)  ,\ J_{\mathbf{s+h}}\left(  \omega\right)  )\ $\ be the
covariance function and let $g_{\parallel\mathbf{h}\parallel}(\omega)$ be of
the form\ given by (\ref{eq 3.7.}). Assume the function $g_{\parallel
\mathbf{h}\parallel}(\omega)$ is characterized by the parameter vector
$\vartheta$. For convenience, we denote the covariance function by
$g_{\parallel\mathbf{h}\parallel}(\omega,\ \vartheta).$ Our object \ here is
to estimate $\underline{\vartheta}$. We note $\ $\ that $\omega$ is the
temporal spectral frequency, and $\left\Vert \mathbf{h}\right\Vert $ is \ the
spatial Euclidean distance. The estimation of the parameters of the covariance
function have also been considered by other authors (see \ for example,
\cite{cressie1999classes}, \cite{gneiting2002nonseparable},
\cite{ma2002spatio}, \cite{ma2003spatio}, \cite{stein2004approximating},
\cite{stein2005statistical}), using either variogram method or likelihood method.

We note that in the case of purely spatial processes, the parameters are
estimated either by minimizing the differences between the estimated
variograms and theoretical variograms evaluated for spatial distances
$\left\Vert \mathbf{h}\right\Vert $ (weighted least squares approach) or by
maximizing the Gaussian likelihood function. Because of the inclusion of
temporal dimension, and if one uses time domain approach, the observations
\ vector to use will be of order $mn\times1,\ $and the variance covariance
matrix of the observation vector will be of dimension $mn\times mn$. The
number of computational operations required \ for inversion of such large
dimensional matrices can be formidable. For example, it is well known that the
calculation of Gaussian likelihood \ from such vectors requires $(nm)^{3}$
operations. In view of this, \cite{stein2004approximating},
\cite{stein2005space} suggested using the restricted likelihood approach, an
extension of the method proposed by \ \cite{Vecchia1988} to reduce the number
of computations. In FV approach proposed here, one does not need inversion
\ of such high dimensional matrices as the likelihood function calculated is
based on complex Gaussianity of Discrete Fourier Transforms evaluated at
several distinct Fourier frequencies and the properties of the DFT's. It is
well known that at these Fourier frequencies, the Discrete Fourier Transforms
of a stationary process are asymptotically independent. Therefore, the
covariance matrix is diagonal. Further, the Discrete Fourier Transforms can be
calculated using the Fast Fourier transform algorithms. It has been shown in
\cite{SubbaRaoa} that the FV estimates are robust against departure from
Gaussianity and are as efficient as Gaussian estimates, if the process happens
to be Gaussian and require less computational time. We \ now define a new
spatio-temporal random process based on differences \ of the observed process
$\left\{  Y_{t}\left(  \mathbf{s}\right)  \right\}  $. Calculate the
differences
\[
X_{ij}\left(  t\right)  =Y_{t}\left(  \mathbf{s}_{i}\right)  -Y_{t}\left(
\mathbf{s}_{j}\right)  ,\quad\ \text{for each \ }t=1,2,\ldots,n,
\]
and for all locations $\mathbf{s}_{i}$, $\mathbf{s}_{j}$ where $\mathbf{s}%
_{i}$ and $\mathbf{s}_{j}$, $(i$ $\neq$\ $j)$ are the pairs that belong to the
set $N(\mathbf{h}_{l})=\left\{  \mathbf{s}_{i},\mathbf{s}_{j}|\ \left\Vert
\mathbf{s}_{i}-\mathbf{s}_{j}\right\Vert =\left\Vert \mathbf{h}_{l}\right\Vert
,\ l=1,2,\ldots,L\right\}  $. Define the Finite Fourier transform of the new
time series $\left\{  X_{ij}\left(  t\right)  |\,\ i\neq j\right\}  $ at the
Fourier frequencies $\omega_{k}=2\pi k/n,\ k=0,\ 1,\ \ldots,\ \left[
n/2\right]  $,
\begin{equation}
J_{X_{ij}}\left(  \omega_{k}\right)  =\frac{1}{\sqrt{(2\pi}n)}\sum
\limits_{t=1}^{n}X_{ij}\left(  t\right)  e^{-it\omega_{k}}=J_{\mathbf{s}_{i}%
}\left(  \omega_{k}\right)  -J_{\mathbf{s}_{j}}\left(  \omega_{k}\right)
\text{, } \label{eq 4.1.}%
\end{equation}
Let $I_{X_{ij}}\left(  \omega_{k}\right)  $ be the second order periodogram of
the time series $\left\{  X_{ij}\left(  t\right)  \right\}  $ given by
\[
I_{X_{ij}}\left(  \omega_{k}\right)  =\left\vert J_{X_{ij}}\left(  \omega
_{k}\right)  \right\vert ^{2}.
\]
\ Let $G_{X_{ij}}(\omega_{k},\ \vartheta)=E(I_{X_{ij}}\left(  \omega
_{k}\right)  )$ be the expectation of $I_{X_{ij}}\left(  \omega_{k}\right)  $.
The function $G_{X_{ij}}(\omega_{k},\ \vartheta)$ is defined as the Frequency
Variogram by \cite{SubbaRaoa}. We can see the similarity of this function to
the classical definition of \ spatio-temporal variogram $2\gamma\left(
\mathbf{h},u\right)  $ defined in section 2 of the present paper (set $u=0$)
in (\ref{eq 1.1.}). The usefulness of FV \ as a measure of dissimilarity
between two spatial processes will be discussed by the authors in a later publication.

From (\ref{eq 4.1.}), we obtain
\begin{align}
E\left[  I_{X_{ij}}\left(  \omega_{k}\right)  \right]   &  =G_{X_{ij}}%
(\omega_{k},\vartheta)\nonumber\\
&  =E\left[  I_{\mathbf{s}_{i}}\left(  \omega_{k}\right)  \right]  +E\left[
I_{\mathbf{s}_{j}}\left(  \omega_{k}\right)  \right]  -2\ Real\ E\left[
I_{\mathbf{s}_{i}\mathbf{s}_{j}}\left(  \omega_{k}\right)  \right]  ,
\label{eq 4.2.}%
\end{align}
where $I_{\mathbf{s}_{i}\mathbf{s}_{j}}\left(  \omega_{k}\right)  $ is the
cross periodogram between the processes $\left\{  Y_{t}\left(  \mathbf{s}%
_{i}\right)  \right\}  $ and $\left\{  Y_{t}\left(  \mathbf{s}_{j}\right)
\right\}  $ and $I_{\mathbf{s}_{i}}\left(  \omega_{k}\right)  $ is the
periodogram of the series $Y_{t}(s_{i})$. For large $n$, it can be shown that
for a stationary process $\ E\left[  I_{\mathbf{s}_{i}}\left(  \omega
_{k}\right)  \right]  =E\left[  I_{\mathbf{s}_{j}}\left(  \omega_{k}\right)
\right]  =g_{0}(\omega_{k};\ \underline{\vartheta})$ and for a stationary and
an isotropic process $E\left[  I_{\mathbf{s}_{i}\mathbf{s}_{j}}\left(
\omega_{k}\right)  \right]  =$ $g_{\left\Vert s_{i}-s_{j}\right\Vert }\left(
\omega_{k};\underline{\vartheta}\right)  $ which is real. Therefore, the
expectation of (\ref{eq 4.2.}) is given by
\begin{equation}
G_{(s_{i},s_{j})}\left(  \omega_{k};\underline{\vartheta}\right)  =2\left[
g_{0}(\omega_{k};\ \underline{\vartheta})-g_{\left\Vert s_{i}-s_{j}\right\Vert
}\left(  \omega_{k};\underline{\vartheta}\right)  \right]  \text{, }
\label{eq 4.3.}%
\end{equation}
It is interesting to compare $G_{(s_{i},s_{j})}\left(  \omega_{k}%
;\underline{\vartheta}\right)  $ with spatio-temporal variogram $\ 2\gamma
\left(  \mathbf{h},u\right)  $ \ given by equation (\ref{eq 1.2.}). The
similarity between these two functions \ shows that one can use the Frequency
variogram which is a frequency domain version of spatio-temporal variogram for
estimating the effective range $||\mathbf{h||,\ }$and also the parameters etc.

Now for the estimation of the parameter vector $\underline{\vartheta}$ we
proceed as in \cite{SubbaRaoa}. Let $M=[\frac{n}{2}]$. Consider the
M-dimensional complex valued random vector,
\[
\underline{\chi}_{\left\Vert s_{i}-s_{j}\right\Vert }\left(  \omega\right)
=\left[  J_{X_{ij}}\left(  \omega_{1}\right)  ,J_{X_{ij}}\left(  \omega
_{2}\right)  ,\ldots,J_{X_{ij}}\left(  \omega_{M}\right)  \right]  ,
\]
which is distributed asymptotically \ as complex normal with mean zero and
with variance covariance matrix with diagonal elements\newline$[g_{\left\Vert
\mathbf{h}\right\Vert }\left(  \omega_{1},\vartheta\right)  ,\ g_{\left\Vert
\mathbf{h}\right\Vert }\left(  \omega_{2},\vartheta\right)  ,\ \ldots
,\ g_{\left\Vert \mathbf{h}\right\Vert }\left(  \omega_{M},\vartheta\right)
],\ $where $||h||=||s_{i}-s_{j}||.$ We note that because of \ asymptotic
independence of Fourier transforms at distinct Fourier frequencies considered
here, the off diagonal elements of the variance covariance matrix of the
complex Gaussian random vector $\underline{\chi}_{\left\Vert \mathbf{h}%
\right\Vert }\left(  \omega\right)  $ are zero. \ Therefore, the minus\ of log
likelihood function can be shown to be proportional to%
\begin{equation}
Q_{n,N\left(  \mathbf{h}\right)  }\left(  \underline{\vartheta}\right)
=\frac{1}{\left\vert N\left(  \mathbf{h}\right)  \right\vert }\sum
\limits_{\left(  \mathbf{s}_{i},\mathbf{s}_{j}\right)  \in N\left(
\mathbf{h}\right)  }\sum\limits_{k=1}^{M}\left[  \ln G_{(\mathbf{s}%
_{i},\mathbf{s}_{j})}\left(  \omega_{k};\underline{\vartheta}\right)
+\frac{I_{Xij}\left(  \omega_{k}\right)  }{G_{(\mathbf{s}_{i},\mathbf{s}_{j}%
)}\left(  \omega_{k};\underline{\vartheta}\right)  }\right]  \text{.}
\label{eq 4.4.}%
\end{equation}
Here $|N\left(  \mathbf{h}\right)  |$ is the total number of all distinct
pairs $\mathbf{s}_{i}$ and $\mathbf{s}_{j}$ such that $N\left(  \mathbf{h}%
\right)  =\left\{  \left(  \mathbf{s}_{i},\mathbf{s}_{j}\right)  |\ \left\Vert
\mathbf{s}_{i}-\mathbf{s}_{j}\right\Vert =||\mathbf{h||}\right\}  $. The above
criterion (\ref{eq 4.4.}) is defined only for one distance $\left\Vert
\mathbf{h}\right\Vert $. Suppose we now define $L$ spatial distances from the
observed data. We can now define an over all criterion for minimization
\begin{equation}
Q_{n}\left(  \vartheta\right)  =\frac{1}{L}\sum\limits_{l=1}^{L}Q_{n,N\left(
\mathbf{h}_{l}\right)  }\left(  \underline{\vartheta}\right)  ,
\label{eq 4.5.}%
\end{equation}
We minimize (\ref{eq 4.5.}) with respect to $\underline{\vartheta}$ (for
details refer to \cite{SubbaRaoa}). The asymptotic normality of the estimator
$\underline{\vartheta}$ obtained by minimizing (\ref{eq 4.5.}) has been proved
in Theorem 2 of the paper of \cite{SubbaRaoa}. To avoid repetition, we refer
to their paper for details. We state the asymptotic distribution of the
estimates. It has been shown in \cite{SubbaRaoa} that under certain
conditions, and for large $n,\ $%
\[
\sqrt{n}\left(  \underline{\vartheta}_{n}-\vartheta_{0}\right)  \overset
{D}{\longrightarrow}N\left(  \underline{0},\nabla^{2}Q_{n}^{-1}\left(
\underline{\vartheta}_{0}\right)  \text{ }V\text{ }\nabla^{2}Q_{n}^{-1}\left(
\underline{\vartheta}_{0}\right)  \right)  ,
\]
where $V=\lim\limits_{n\rightarrow\infty}var\left[  \frac{1}{\sqrt{n}}\nabla
Q_{n}\left(  \vartheta_{0}\right)  \right]  $, $\nabla Q_{n}\left(
\vartheta_{0}\right)  $ is a vector of first order partial derivatives,
$\nabla^{2}Q_{n}\left(  \vartheta_{0}\right)  $ is the matrix of second order
partial derivatives. We consider the estimation, prediction etc. for a real
example in the following section.

\section{Real Data Analysis}

For our illustration, we consider the Air Pollution data analyzed by
\cite{sahu2005bayesian}, and we refer to their paper for full details. The
data analyzed corresponds to atmospheric \ particulate matter that is less
than $2.5$ $\omega m$ in size (usually known as PM$_{2.5}$) \ which is one of
six primary air pollutants and is a mixture of fine particles and gaseous
compounds \ such as sulphur dioxide (SO$_{2}$) and nitrogen oxides. The data
was observed at $15$ monitoring stations \ in New York city \ during the first
$9$ months of the year 2002. The data was observed once in every $3$ days,
thus giving $91$ equally spaced time series for each monitoring station. The
total number of observations are $1365=15\times91$. The data can be obtained
from the website \ \ \ http://www.blackwellpublishing.com/rss. We use the data
given at the $15$ locations along with their spatial coordinates. The spatial
coordinates of $625$ nearby locations are also known and can be found in the
website. We also consider the estimation of the data at these locations.
\begin{figure}
[ptb]
\begin{center}
\includegraphics[
natheight=9.821700in,
natwidth=8.523600in,
height=5.9698in,
width=5.1828in
]%
{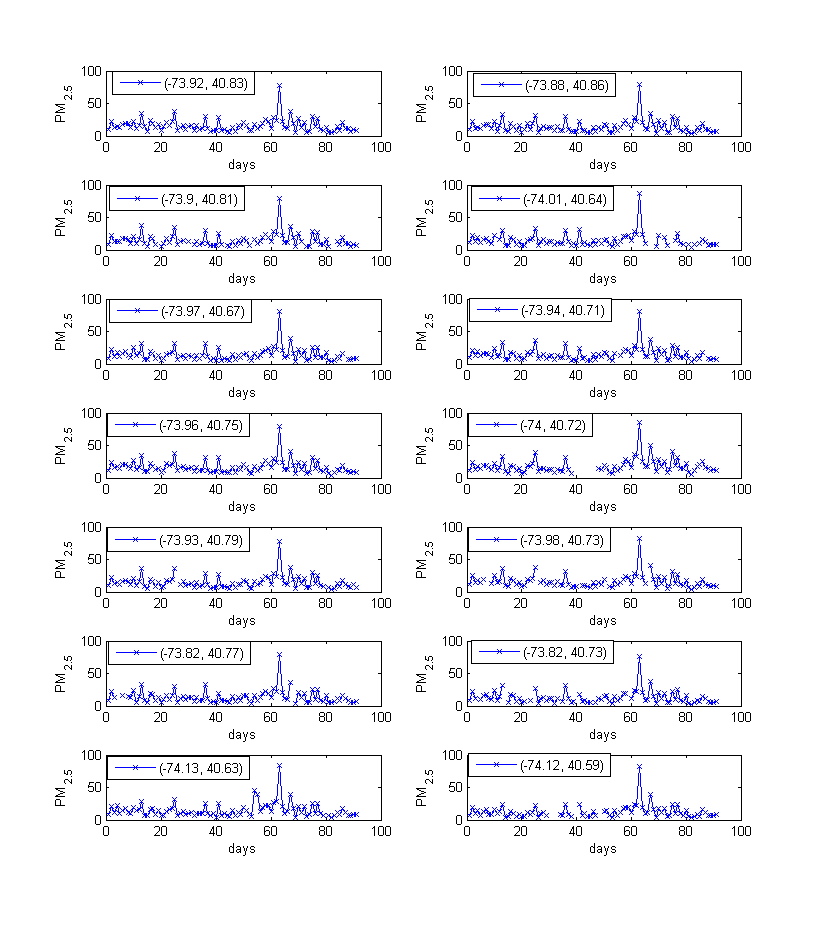}%
\caption{Time Series Data at the 14 locations with their spatial coordinates}%
\label{14ser}%
\end{center}
\end{figure}

Out of $1365$ data points, $126$ were missing \ and the missing values are
estimated as follows. Let $X\left(  \mathbf{s}_{i},t\right)  $
\ $(i=1,\ 2,\ \ldots15;$ $t=1,\ 2,\ \ldots91)$ denote the PM$_{2.5}$
\ observation at the location $\mathbf{s}_{i}$ and at time $t$ and suppose the
observation $X\left(  \mathbf{s}_{0},t_{0}\right)  $ is missing. The missing
observation is estimated by averaging over the data at other locations at time
$t_{0}$ by $\widetilde{X}\left(  \mathbf{s}_{0},t_{0}\right)  =(1/14)\sum
_{j=1,\neq0}^{14}X\left(  \mathbf{s}_{_{j}},t_{0}\right)  $.

We used the first differences $\{Y_{t}\left(  \mathbf{s}\right)  =\Delta
_{t}X\left(  \mathbf{s},t\right)  :\mathbf{s}\in\mathbb{R}^{2},\ \ t\in
\mathbb{Z\}}$ to remove the linear trend as suggested by
\cite{sahu2005bayesian}, and used the detrended data for our analysis. The
time series data from the $14$ locations are plotted in Figure \ref{14ser}.
The location $11$ (with coordinates: latitude $-73.84$, longitude $40.77$) has
large number of missing values (only the first 23 observations are available)
and, therefore, we have chosen to estimate all 91 observations at this
location using the data from other 14 locations. We compared the estimated
values with the 23 observed .\ The sample auto-correlation plots of the $14$
series have shown strong correlation at lag $1$. On the basis of BIC
criterion, we decided that the $14$ time series can adequately be modelled by
an AR(1) model. The corresponding second order spectrum is \ \ $g_{0}\left(
\omega,\vartheta\right)  =\sigma^{2}{\large /}\left\vert 1+\varphi
z\right\vert ^{2}$ where $z=\exp(i\omega)$ and the vector of the parameters
$\ $are$\ \vartheta$ $=(\sigma^{2},\phi)$ \ and these are estimated \ by
minimizing the criterion (\ref{eq 4.5.}) with $L=91$. We note that we scaled
the equation (\ref{eq 3.4.}) such that $\sigma_{e}^{2}=\sigma^{2}$, see
(\ref{eq 3.7.0}).

The final estimates \ of the AR model have been found to be $\widehat{\varphi
}=0.4659$; \ $\widehat{\sigma}=5.9264$.\ Using these estimated values, all the
elements of the vector $G_{0}\left(  \omega\right)  $ and the elements of
$\ $the square matrix $F_{m}\left(  \omega\right)  $ (which is of order
$14\times14$) are evaluated. The vectors $J_{0}(\omega)$ at the \ Fourier
frequencies \ $\omega_{k}=2\pi k/2^{6}$\ are estimated using the equation
$\underline{\widehat{G}}_{0}^{\prime}\left(  \omega_{k}\right)  \widehat
{F}_{m}^{-1}\left(  \omega_{k}\right)  \underline{J}_{m}\left(  \omega
_{k}\right)  $. The data at the location $11$ are estimated using the equation
(\ref{eq 5.1. a}). The plot of the 91 estimated values (with (+) sign), plot
of the first $23$ given observations (with $o$ sign), corresponding $95\%$
confidence bands using (\ref{VarPred}) are given in Fig 2.\ We see a good
agreement between the estimated values and the observed values, suggesting
that \ the prediction methodology works very well in this case. Also we find
that the strong spatial correlation and the temporal correlation \ can
satisfactorily be described by the spatio-temporal covariance function
\ defined here.
\begin{figure}
[ptbh]
\begin{center}
\includegraphics[
natheight=4.708000in,
natwidth=6.520700in,
height=3.627in,
width=5.0142in
]%
{Pred_CIn11a.png}%
\caption{Location 11. Plot of the observed values (denoted by + sign),
predicted values (denoted by 'o' sign), 95\% prediction intervals (denoted by
-- sign)}%
\label{Pred11}%
\end{center}
\end{figure}

In order to check \ the overall performance, we computed the leave-one-out
cross-validation (\cite{giraldo2010continuous}) criterion . Here we estimated
the data at one location, taken one at a time, using the data given at other
13 locations. The Mean Square Error calculated for all the 14 locations is
\[
MSSE=\sum SSE\left(  j\right)  /(14\ast91)=15.3050,
\]
where
\[
SSE\left(  j\right)  =\sum_{t=1}^{91}\left(  Y_{t}\left(  \mathbf{s}%
_{j}\right)  -\widehat{Y}_{t}\left(  \mathbf{s}_{j}|\mathbf{s}_{k}%
\neq\mathbf{s}_{j}\right)  \right)  ^{2},
\]
and $\widehat{Y}_{t}\left(  \mathbf{s}_{j}|\mathbf{s}_{k}\neq\mathbf{s}%
_{j}\right)  $ is the estimator of the data at time $t$ at location
$\mathbf{s}_{j}$ conditional on the data at the locations $\left\{
\mathbf{s}_{1},\mathbf{s}_{2},\ldots\mathbf{s}_{j-1},\mathbf{s}_{j+1}%
,\ldots\mathbf{s}_{13}\right\}  $. We have also estimated the prediction error
variance using the equation (21) \ for each location, and they are shown in
Fig 3 (the larger the diameter of the circle, the higher the variance).

We estimated the PM$_{2.5}$ values for all $t=1,2,\ldots91$ and for all the
$625$ locations (including $14$ locations where the data is available). The
plot of the these predicted values are given in Fig 4..
\begin{figure}
[ptb]
\begin{center}
\includegraphics[
natheight=5.354000in,
natwidth=6.801800in,
height=3.7308in,
width=4.734in
]%
{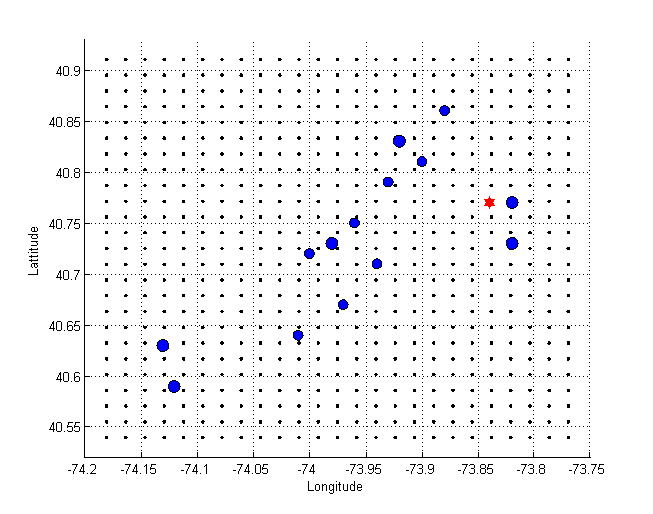}%
\caption{Grid of the 625 locations. Dark circles correspond to the locations
(data available).}%
\label{Loc625}%
\end{center}
\end{figure}
%

\begin{figure}
[ptb]
\begin{center}
\includegraphics[
natheight=5.354000in,
natwidth=6.801800in,
height=3.96in,
width=5.0263in
]%
{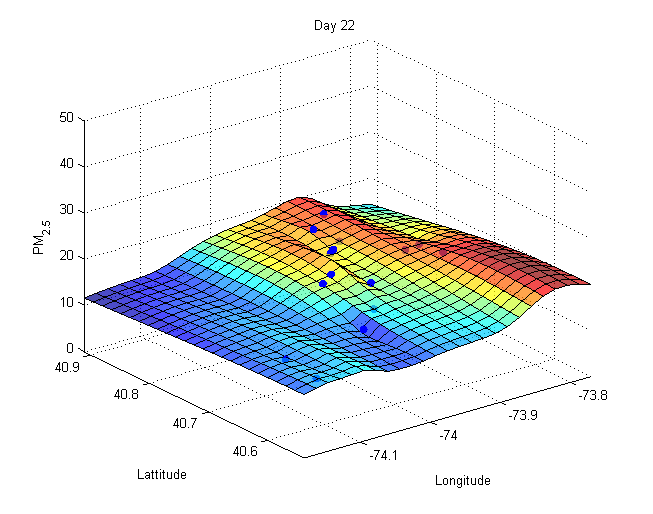}%
\caption{Estimated PM$_{2.5}$ for all the locations at $t=22$. }%
\label{Day22}%
\end{center}
\end{figure}

To check the forecasting performance of the method described \ here we
considered the prediction at the location \#8. We considered the data,
$t=1,2,...83$ as given and estimated the rest of the values ($84,83,\ldots91$)
using the AR(1) model fitted. The forecasting methodology is well known, and
therefore, we briefly summarise the method (details can be found in any time
series book).

Let%
\[
\widehat{Y}_{n+h}\left(  \mathbf{s}_{0}|t=1,2,\ldots,n;\mathbf{s}%
_{j};j=1,2,\ldots m\right)  =E\left(  Y_{n+h}\left(  \mathbf{s}_{0}\right)
|Y_{t}\left(  \mathbf{s}_{j}\right)  ;t=1,2,\ldots n,\mathbf{s}_{j}%
;j=1,2,\ldots m\right)  \ \
\]
For the AR(1) model fitted, we have
\[
\widehat{Y}_{n+h}\left(  \mathbf{s}_{0}|t=1,2,\ldots n,\mathbf{s}%
_{j};j=1,2,\ldots m\right)  =\varphi^{h}\widehat{Y}_{n}\left(  \mathbf{s}%
_{0}\right)  ,
\]
hence
\[
Y_{n+h}\left(  \mathbf{s}_{0}\right)  -\widehat{Y}_{n+h}\left(  \mathbf{s}%
_{0}\right)  =\varphi^{h}\left(  Y_{n+h}\left(  \mathbf{s}_{0}\right)
-\widehat{Y}_{n}\left(  \mathbf{s}_{0}\right)  \right)  +\sum_{k=1}^{h}%
\varphi^{k-1}\varepsilon_{n+k}\left(  \mathbf{s}_{0}\right)  ,
\]
where $\widehat{Y}_{n}\left(  \mathbf{s}_{0}\right)  $ denotes the spatial
prediction and $\varepsilon_{t}\left(  \mathbf{s}_{0}\right)  $ is the
innovation error at $\mathbf{s}_{0}$ and at time $t$. Therefore, we obtain%
\begin{align}
E\left(  Y_{n+h}\left(  \mathbf{s}_{0}\right)  -\widehat{Y}_{n+h}\left(
\mathbf{s}_{0}\right)  \right)  ^{2} &  =\varphi^{2h}E\left(  Y_{n}\left(
\mathbf{s}_{0}\right)  -\widehat{Y}_{n}\left(  \mathbf{s}_{0}\right)  \right)
^{2}+\sigma^{2}\frac{\varphi^{2h}-1}{\varphi^{2}-1}\nonumber\\
&  =\varphi^{2h}\int\limits_{-\pi}^{\pi}\sigma_{m}^{2}\left(  \omega\right)
d\omega+\sigma^{2}\frac{\varphi^{2h}-1}{\varphi^{2}-1}.\label{VarPred}%
\end{align}%
\begin{figure}
[ptb]
\begin{center}
\includegraphics[
natheight=4.343900in,
natwidth=5.593600in,
height=3.6443in,
width=4.6864in
]%
{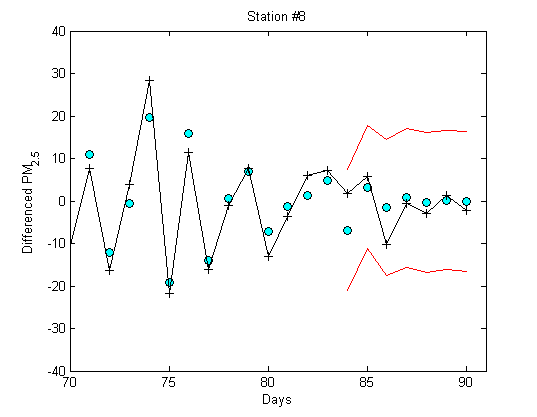}%
\caption{Location 8. Predicted values ('o'). observed values (+), Prediction
intervals (--) }%
\label{FigPredSt8}%
\end{center}
\end{figure}

\section{{\protect\Large Appendix: }Discrete Fourier Transforms and their
properties}

Let us assume that we have time series data from $m$ locations spatially
distributed. \newline Let $\left\{  Y_{t}\left(  \mathbf{s}_{i}\right)
|i=1,2,\ldots m;\ t=1,\ldots,n\right\}  $, be a sample from the zero mean
spatio-temporal stationary process $\left\{  Y_{t}\left(  \mathbf{s}\right)
\right\}  $. \ Consider the time series data $\left\{  Y_{t}\left(
\mathbf{s}_{i}\right)  |t=1,\ldots,n\right\}  $ at the location $\mathbf{s}%
_{i}$, and define the Discrete Fourier transform (DFT)
\begin{equation}
J_{\mathbf{s}_{i}}\left(  \omega_{k}\right)  =\frac{1}{\sqrt{2\pi n}}%
\sum\limits_{t=1}^{n}Y_{t}\left(  \mathbf{s}_{i}\right)  e^{-it\omega_{k}},
\label{eq 2.1.}%
\end{equation}
where $\omega_{k}=\frac{2\pi k}{n},\ \ k=0,\ \pm1,\ \ldots,\ \pm\left[
\frac{n}{2}\right]  $. In practice one uses Fast Fourier Transform algorithm
to compute the DFT. \ It is well known that the number of operations required
to calculate $n$ discrete Fourier transforms from $n$ dimensional data is of
order $\ n\log n$. The corresponding second order periodogram \ is defined as
\[
I_{\mathbf{s}_{i}}\left(  \omega_{k}\right)  =\left\vert J_{\mathbf{s}_{i}%
}\left(  \omega_{k}\right)  \right\vert ^{2},
\]
and the cross periodogram between the time series $\left\{  Y_{t}\left(
\mathbf{s}_{i}\right)  \right\}  $ \ and $\left\{  Y_{t}\left(  \mathbf{s}%
_{j}\right)  \right\}  $ is given by $I_{\mathbf{s}_{i}s_{j}}\left(
\omega_{k}\right)  =J_{\mathbf{s}_{i}}\left(  \omega_{k}\right)
J_{\mathbf{s}_{j}}^{\ast}\left(  \omega_{k}\right)  $.

It is well known that the periodogram \ is an asymptotically unbiased
estimator of the second order spectral density function, but it is not mean
square consistent and hence \ to obtain a consistent estimate the periodograms
are smoothed using various kernels (see \cite{Priestley_book-81}).\ It is well
known \ that (see \ for example \cite{Priestley_book-81})%
\begin{align}
E\left(  J_{\mathbf{s}_{i}}\left(  \omega\right)  \right)   &  =0,\nonumber\\
Var\left(  J_{\mathbf{s}_{i}}\left(  \omega\right)  \right)   &  =E\left(
I_{\mathbf{s}_{i}}\left(  \omega\right)  \right)  \simeq g_{\mathbf{s}_{i}%
}\left(  \omega\right)  =g_{0}(\omega), \label{eq 2.2.}%
\end{align}
where $g_{0}(\omega)$ is the second order spectral density function of the
random process.\ We further note\ that $\sigma_{Y}^{2}=\int g_{0}%
(\omega)d\omega$. From (\ref{eq 2.1.}), by inversion we get%
\begin{equation}
Y_{t}\left(  \mathbf{s}\right)  =\sqrt{\frac{n}{2\pi}}\int\limits_{-\pi}^{\pi
}J_{\mathbf{s}}\left(  \omega\right)  e^{it\omega}d\omega\text{.}
\label{Y_invF}%
\end{equation}
The above relation (given by the (\ref{Y_invF})) shows that the process
$Y_{t}\left(  \mathbf{s}\right)  $\ can be decomposed into sine and cosine
functions, and the amplitude $J_{\mathbf{s}}\left(  \omega\right)  $
measures\ the contribution made by the harmonic component with frequency
$\omega$ to the estimated total power $\sigma_{Y}^{2}$. It \ is also well
known that the power spectrum $g_{0}(\omega)$ \ is a frequency decomposition
of the power $\sigma_{Y}^{2}$ of the signal \ $\{Y_{t}\left(  \mathbf{s}%
_{i}\right)  \}$ and, therefore, performing spectral analysis on the data is
considered to be equivalent to the classical analysis of variance (ANOVA),
where one checks for the contribution of each component in the regression
model via ANOVA decomposition.

In the following proposition we summarise the well known properties of the
Fourier Transforms of stationary process. For details refer to
\cite{Brill-book-01}, \cite{Priestley_book-81}, \cite{dwivedi2011test}.

\begin{proposition}
Let $J_{\mathbf{s}_{i}}\left(  \omega_{k}\right)  $ and $J_{\mathbf{s}_{j}%
}\left(  \omega_{k}\right)  $ be the discrete Fourier Transforms of the
spatio-temporal stationary processes $\{Y_{t}\left(  \mathbf{s}_{i}\right)
\},\ \{Y_{t}\left(  \mathbf{s}_{j}\right)  \}$ respectively. For large $n$,
\[
Cov\left(  J_{\mathbf{s}_{i}}\left(  \omega_{k}\right)  ,J_{\mathbf{s}_{i}%
}\left(  \omega_{k^{^{\prime}}}\right)  \right)  \simeq0,\quad k\neq
k^{\prime},
\]%
\begin{align}
Cov\left(  J_{\mathbf{s}_{i}}\left(  \omega_{k}\right)  ,J_{\mathbf{s}_{j}%
}\left(  \omega_{k}\right)  \right)   &  =E\left[  I_{\mathbf{s}_{i}%
\mathbf{s}_{j}}\left(  \omega_{k}\right)  \right] \nonumber\\
&  \simeq\frac{1}{2\pi}\sum\limits_{n=-\infty}^{\infty}c\left(  \mathbf{s}%
_{i}-\mathbf{s}_{j},n\right)  e^{-in\omega_{k}}=g_{\mathbf{s}_{i}%
-\mathbf{s}_{j}}\left(  \omega_{k}\right)  , \label{eq 2.3.}%
\end{align}
where $g_{\left(  \mathbf{s}_{i}-\mathbf{s}_{j}\right)  }\left(
\omega\right)  $ is defined as the cross spectrum between the two processes
and it is usually a complex valued function. If the process is isotropic,
then
\[
c\left(  \mathbf{s}_{i}-\mathbf{s}_{j},n\right)  =c\left(  \left\Vert
\mathbf{s}_{i}-\mathbf{s}_{j}\right\Vert ,n\right)  =c\left(  \left\Vert
\mathbf{s}_{i}-\mathbf{s}_{j}\right\Vert ,-n\right)  \text{.}%
\]
and under the isotropy assumption the cross spectrum $g_{\mathbf{s}%
_{i}-\mathbf{s}_{j}}\left(  \omega\right)  $ between the two processes reduces
to
\[
g_{\parallel\mathbf{h}\parallel}(\omega)=\frac{1}{2\pi}\sum\limits_{n=-\infty
}^{\infty}c(\left\Vert \mathbf{h}\right\Vert ,n)e^{-in\omega},\quad\left\vert
\omega\right\vert \leq\pi,
\]
and the spectral function $g_{\parallel\mathbf{h}\parallel}(\omega)$, where
$\ \left\Vert \mathbf{h}\right\Vert =\left\Vert \mathbf{s}_{i}-\mathbf{s}%
_{j}\right\Vert ,\ $is \ real and symmetric in $\omega$.
\end{proposition}

\begin{proof}
The above results are well known and hence details are omitted.
\end{proof}

If the random process $\left\{  Y_{t}\left(  \mathbf{s}_{i}\right)  \right\}
$ is Gaussian, then the complex valued random variables \newline$\left\{
J_{\mathbf{s}_{i}}\left(  \omega_{k}\right)  |\ k=0,1,\ldots,\left[  \frac
{n}{2}\right]  \right\}  $ will be asymptotically independent, \ and will be
distributed as complex Gaussian, \ and each $J_{\mathbf{s}_{i}}\left(
\omega_{k}\right)  $ will be distributed as complex normal with mean zero and
variance proportional to $g_{0}\left(  \omega_{k}\right)  $.

As pointed out earlier the second order periodogram $\ I_{\mathbf{s}_{i}%
}(\omega)$ defined above is always real, whereas the cross periodogram
$I_{\mathbf{s}_{i},\mathbf{s}_{j}}\left(  \omega_{k}\right)  $ defined above
between the spatial processes $\left\{  Y_{t}\left(  \mathbf{s}_{i}\right)
\right\}  $ and $\left\{  Y_{t}\left(  \mathbf{s}_{j}\right)  \right\}  $ \ is
usually a complex valued function. Under the isotropy assumption, however, the
cross spectrum is a \ function of the Euclidean distance $\ \left\Vert
\mathbf{h}\right\Vert =\left\Vert \mathbf{s}_{i}-\mathbf{s}_{j}\right\Vert $
and the temporal frequency $\omega$ (i.e., spectral in time, but not in space)
and, therefore, it is real. It is interesting to see the similarity between
the above function $\ $and the spectral density functions defined by
\cite{cressie1999classes} and \cite{stein2005space} which are spectral in time
but not in space.

In the following proposition we show that the Discrete Fourier Transform \ of
$Y_{t}\left(  \mathbf{s}\right)  $ can be written in terms of \ orthogonal set
function $Z_{Y}\left(  \underline{\lambda},\omega\right)  $.

\begin{proposition}
Let $\ J_{\mathbf{s}}\left(  \omega\right)  $ be the Discrete Fourier
Transform of $\{Y_{t}(s)\}$ and let the spectral representation of $Y_{t}(s)$
be given by (3). Then%
\begin{equation}
J_{\mathbf{s}}\left(  \omega\right)  \simeq\int e^{i\underline{s}%
\,\underline{\lambda}}\sqrt{\frac{n}{2\pi}}dZ_{Y}\left(  \underline{\lambda
},\omega\right)  . \label{eq 2.8.}%
\end{equation}
\bigskip
\end{proposition}

\begin{proof}
Substitute the spectral representation (3) for $Y_{t}\left(  \mathbf{s}%
\right)  $ in (\ref{eq 2.1.}), and after some simplification, we obtain
\begin{equation}
J_{\mathbf{s}}\left(  \omega\right)  =\int\int e^{i\underline{s}%
\,\underline{\lambda}}\left[  e^{i\left(  n+1\right)  \frac{\varphi}{2}}%
F_{n}^{\frac{1}{2}}\left(  \varphi\right)  \right]  dZ_{Y}\left(
\underline{\lambda},\omega\right)  , \label{eq 2.7.}%
\end{equation}
where $\varphi=\mu-\omega$, $\int$ is a $d$ dimensional multiple integral,
(see \cite{Priestley_book-81}, p. 419) and in obtaining the above, we used the
result%
\[
\sum\limits_{t=1}^{n}e^{it\varphi}=e^{i\left(  n+1\right)  \frac{\varphi}{2}%
}\left[  \frac{\sin n\frac{\varphi}{2}}{\sin\frac{\varphi}{2}}\right]
=e^{i\left(  n+1\right)  \frac{\varphi}{2}}\sqrt{2\pi n}F_{n}^{\frac{1}{2}%
}\left(  \varphi\right)  \text{, }%
\]
where the Fej\'{e}r kernel $F_{n}\left(  \varphi\right)  $ is given by%
\[
F_{n}\left(  \varphi\right)  =\frac{1}{2\pi n}\frac{\sin^{2}n\frac{\varphi}%
{2}}{\sin^{2}\frac{\varphi}{2}}\text{.}%
\]
It is well known that the Fej\'{e}r kernel behaves like a Dirac Delta function
as $n\rightarrow\infty$ and as $\varphi\rightarrow0$, $F_{n}\left(
\varphi\right)  =\mathit{O}\left(  n\right)  $. As pointed out by
\cite{Priestley_book-81}, p. 419), that $F_{n}^{\frac{1}{2}}\left(
\varphi\right)  $ does not strictly tend to a Dirac Delta $\delta-$function as
$n\rightarrow\infty$, nevertheless, behaves in a similar manner to a $\delta
-$function. In particular as $n\rightarrow\infty$ and for all $\varphi\neq0$,
$F_{n}^{\frac{1}{2}}\left(  \varphi\right)  \rightarrow0$, and as
$\varphi\rightarrow0$, $F_{n}^{\frac{1}{2}}\left(  \varphi\right)
\rightarrow\sqrt{n/2\pi}$. Therefore, as $n\rightarrow\infty$, $F_{n}%
^{\frac{1}{2}}\left(  \varphi\right)  $ vanishes everywhere except at the
origin. In view of this, for large $n$, we \ have the result%
\[
J_{\mathbf{s}}\left(  \omega\right)  \simeq\int e^{i\underline{s}%
\,\underline{\lambda}}\sqrt{\frac{n}{2\pi}}dZ_{Y}\left(  \underline{\lambda
},\omega\right)  .
\]
We note\ that the above integral is over the wave number space $\lambda$ only.
\end{proof}

\begin{proposition}
Let $\left\{  e_{t}(\mathbf{s)}|\mathbf{s}\in\mathbb{R}^{d},t\in
\mathbb{Z}\right\}  $ be a white noise process in space and time, that is, we
assume the process has constant spectrum, it satisfies the following
conditions.%
\begin{align*}
E\left(  e_{t}\left(  \mathbf{s)}\right)  \right)   &  =0,\\
Var(e_{t}\left(  \mathbf{s)}\right)   &  =\sigma_{e}^{2}\text{, does not
depend on }\mathbf{s}\text{ or }t,\\
Cov\left(  e_{t}\left(  \mathbf{s)}\right)  ,e_{t^{^{\prime}}}\left(
\mathbf{s)}\right)  \right)   &  =\sigma_{e}^{2}I\left(  \mathbf{s}%
,\mathbf{s}^{\prime}\right)  I\left(  t,t^{\prime}\right)  ,
\end{align*}
where%
\begin{align*}
I(\mathbf{s},\mathbf{s}^{\prime})  &  =\left\{
\begin{array}
[c]{c}%
1\text{ if }\mathbf{s}=\mathbf{s}^{\prime},\\
0\text{ otherwise, }%
\end{array}
\right. \\
I(t,t^{\prime})  &  =\left\{
\begin{array}
[c]{c}%
1\text{ if }t=t^{\prime},\\
0\text{ otherwise.}%
\end{array}
\right.
\end{align*}
Let the Discrete Fourier transform of the white noise process be%
\[
J_{\mathbf{s},e}\left(  \omega\right)  =\frac{1}{\sqrt{2\pi n}}\sum
\limits_{t=1}^{n}e_{t}(\mathbf{s)}e^{-it\omega}.
\]
Then%
\begin{equation}
J_{\mathbf{s},e}\left(  \omega\right)  \simeq\int e^{i\mathbf{s}%
\cdot\underline{\lambda}}\left[  \sqrt{\frac{n}{2\pi}}\right]  dZ_{e}\left(
\underline{\lambda},\omega\right)  , \label{DFT_WN}%
\end{equation}
where the orthogonal random process $Z_{e}\left(  \underline{\lambda}%
,\omega\right)  $ satisfies
\begin{align*}
EdZ_{e}\left(  \underline{\lambda},\omega\right)   &  =0,\\
E\left\vert dZ_{e}\left(  \underline{\lambda},\omega\right)  \right\vert ^{2}
&  =\frac{\sigma_{e}^{2}}{\left(  2\pi\right)  ^{d+1}}d\underline{\lambda
}d\omega\text{.}%
\end{align*}

\end{proposition}

\begin{proof}
It is similar to Proposition 2 and hence omitted.
\end{proof}

\bigskip

\textbf{Acknowledgement.} \textit{The publication was supported by the
T\'{A}MOP-4.2.2.C-11/1/KONV-2012-0001 project. The project has been supported
by the European Union, co-financed by the European Social Fund. The visit of
Subba Rao to the CRRAO AIMSCS was supported by a grant from the Department of
Science and Technology, Government of India, grant number SR/S4/516/07. The
authors are thankful to the editor, co-editor and the referees for their
valuable comments and suggestions. They are also thankful to Dr Suhasini Subba
Rao, Texas A\&M University, \ USA for her suggestions for improvement of the
paper.}

\bibliographystyle{plain}
\bibliography{00BiblMM13}

\end{document}